\documentclass[11pt]{article}
\usepackage{amssymb}
\usepackage{amsthm}
\usepackage{amsmath,amsfonts,amssymb}

\usepackage[mathscr]{eucal}
\usepackage{exscale}
\usepackage{natbib}
\usepackage{bm}
\usepackage{eqlist}
\usepackage[dvipsnames]{color}
\usepackage[Lenny]{fncychap}
\usepackage{multirow}
\usepackage{graphicx}
\usepackage{algpseudocode}
\usepackage{algorithm}
\usepackage{caption}
\usepackage{hyperref}
\usepackage{framed}
\usepackage{color}
\usepackage{booktabs}

\usepackage{framed}
\usepackage{algpseudocode}

\usepackage{tabu}
\usepackage{tikz}

\usepackage[utf8]{inputenc} 
\usepackage[english]{babel} 
\usepackage[T1]{fontenc}    
\usepackage{amssymb,amsmath,amsfonts} 
\usepackage{braket} 
\usepackage{caption}
\usepackage{algorithm}

\hypersetup{
    colorlinks=true,
    linkcolor=blue,
    citecolor=blue,
    citebordercolor={1 1 0},
}

\newtheorem{theorem}{Theorem}[section]

\newtheorem{definition}[theorem]{Definition}
\newtheorem{remark}[theorem]{Remark}

\font\bigbf=cmbx10 scaled \magstep3

\begin{document}

\title{\bigbf  Exploring the sensing power of mixed vehicle fleets}

\author{Ke Han$^{a}\thanks{Corresponding author, e-mail: kehan@swjtu.edu.cn;}$
\quad Wen Ji$^{b}$\quad Yu (Marco) Nie$^{c}$\quad Zhexian Li$^{d}$\quad Shenglin Liu$^{b}$\\\\
 $^{a}$\textit{\small School of Economics and Management,}\\
\textit{\small Southwest Jiaotong University, Chengdu, Sichuan 610031, China}\\
$^{b}$\textit{\small School of Transportation and Logistics,}\\
\textit{\small Southwest Jiaotong University, Chengdu, Sichuan 611756, China}\\
$^{c}$\textit{\small Department of Civil and Environmental Engineering,}\\
\textit{\small Northwestern University, Evanston, IL 60208-3109, USA}\\
$^{d}$\textit{\small Department of Civil and Environmental Engineering,}\\
\textit{\small University of Southern California, Los Angeles, CA 90089, USA}
}

\maketitle

\begin{abstract}
Vehicle-based mobile sensing, also known as drive-by sensing, efficiently surveys urban environments at low costs by leveraging the mobility of urban vehicles. While recent studies have focused on drive-by sensing for fleets of a single type, our work explores the sensing power and cost-effectiveness of a mixed fleet that consists of vehicles with distinct and complementary mobility patterns. We formulate the drive-by sensing coverage (DSC) problem, proposing a method to quantify sensing utility and an optimization procedure that determines fleet composition, sensor allocation, and vehicle routing for a given budget. Our air quality sensing case study in Longquanyi District (Chengdu, China) demonstrates that using a mixed fleet enhances sensing utilities and achieves close approximations to the target sensing distribution at a lower cost. Generalizing these insights to two additional real-world networks, our regression analysis uncovers key factors influencing the sensing power of mixed fleets. This research provides quantitative and managerial insights into drive-by sensing, showcasing a positive externality of urban transport activities.
\end{abstract}


\section{Introduction}\label{secIntro}

In the era of smart cities, the advent of ubiquitous sensing devices (e.g. smart phones, low-cost wireless sensors) has given rise to a data collection paradigm known as {\it mobile crowdsensing} \citep{DSXVF2019}, which has offered unprecedented opportunities to boost sensing capabilities by leveraging human mobility \citep{MZY2014, LBDNT2005}. Vehicle-based mobile crowdsensing, also known as drive-by sensing (DS), has gained popularity because it is relatively cheap and the sensors can be easily moved around. 

DS has seen a wide range of smart city applications such as air quality monitoring \citep{Kumar2015, SHS2021}, traffic state estimation \citep{LSLHLW2009, ZLZLZ2013, DCYLGS2015}, noise pollution monitoring \citep{Alsina-Pages2017, Cruz2020} and infrastructure health inspection \citep{Eriksson2008, Wang2014}. Its sensing quality depends on the spatial-temporal coverage, typically measured by the number of recorded observations in an area during a certain period of time \citep{JHL2023}. The spatial-temporal coverage and resolution, as well as flexibility and reliability of DS, are closely related to the mobility patterns of the vehicular fleet used for sensing. Each fleet type has distinct set-up, maintenance and management costs.  Taxis and buses are the most frequently considered in the literature  \citep{Meegahapola2019, Kaivonen2020, JHL2023}, among other choices including dedicated vehicles\footnote{Dedicated vehicles can be freely scheduled and navigated to perform sensing tasks.} \citep{FZGJGW2021, JHG2023} and commercial trucks \citep{ADRMdR2018, dADKKR2020}.

\begin{figure}[h!]
\centering
\includegraphics[width=\textwidth]{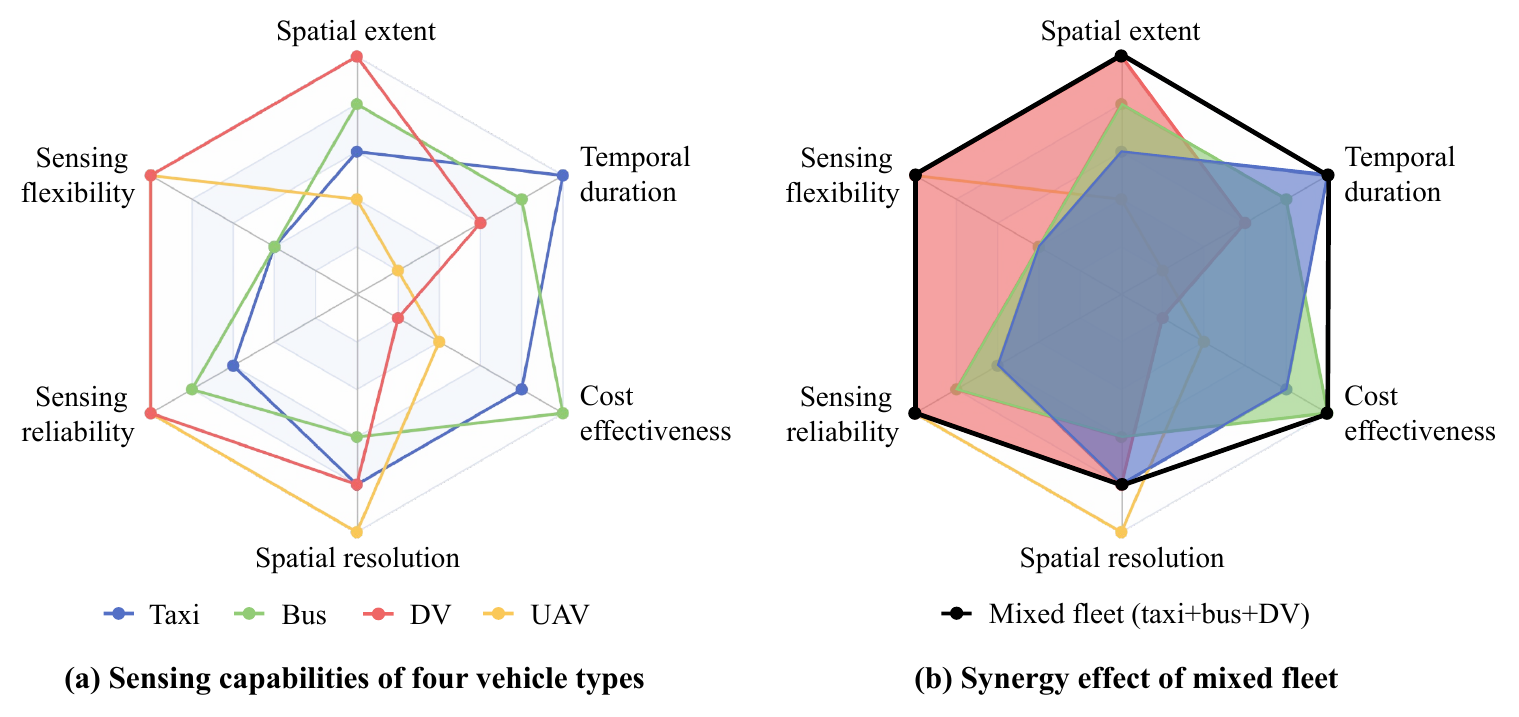}
\caption{Characterization of drive-by sensing capabilities of different vehicle types. Spatial extent: The size of area within the reach of the fleet; Temporal duration: Maximum time span for continuous scanning; Sensing reliability: Likelihood of a given area being scanned at least once within certain period; Spatial resolution: Smallest spatial unit that can be scanned; Sensing flexibility: Maneuverability of the fleet to cover a target area; Cost effectiveness: Procurement, operational and maintenance costs to collect a unit amount of data.}
\label{figthreetypes}
\end{figure}

Existing DS studies have attempted to quantify sensing power via empirical analyses \citep{OASSR2019}, optimize sensing utility through sensor allocation \citep{CLGZX2017, DH2023, JHL2023partC}, or devise incentives to improve vehicle distribution \citep{XCPJZN2019}. The vast majority is focused on the use of a fleet with the same type of vehicles. However, the sensing capability of a homogenous fleet has limitations, such as unbalanced spatial distributions, limited operating hours, and high operational costs. Figure \ref{figthreetypes}(a) \citep{JHL2023} compares the six aspects of sensing capability of different vehicle types. In this paper, by integrating fleets taxis, buses and dedicated vehicles (DVs), we explore the synergy between their sensing capabilities.  These three fleet types are chosen as they are widely available in metropolitan areas, and their capabilities seem complementary to each other, as shown in Figure \ref{figthreetypes}(b). 

Specifically, given a fixed budget, we aim to determine the number of sensors, fleet composition, and operational parameters to achieve optimal spatial-temporal coverage, subject to various mobility and operational constraints. This is coined the {\it drive-by sensing coverage} (DSC) problem, which builds on two components: 
\begin{enumerate}
\item[(1)] {\bf Optimization objective}: The sensing efficacy is quantified in space and time using a general utility function, which accommodates various sensing requirements and user-defined preferences. Notably, an ideal sensing distribution can be derived from such an objective via first-order analysis. 

\item[(2)] {\bf Mobility constraints}: The maximization of sensing efficacy is constrained by the spatial-temporal distribution of the sensing fleet, which is the superposition of  heterogenous mobility patterns of the three vehicle types. 
\end{enumerate}

The main challenge is that the movements of taxis, buses and DVs are associated with distinct levels of uncertainty and controllability, which has two implications: (1) It is difficult to express the mobility constraints analytically, and a heuristic approach is needed to solve for the three vehicle types simultaneously; (2) Besides the number of instrumented vehicles in each type, the problem needs to determine additional operational parameters, namely number of sensors per bus line, and vehicle routes for the DVs.  To tackle this challenge, we design an optimization procedure that alternates between two coupled subproblems: (i) A taxi-bus subproblem, which is formulated using a probabilistic and time-averaged approach and solved as a convex program; and (ii) A DV subproblem, which is formulated as a team orienteering problem with accumulable utilities, and solved using a dual-spatial-scale greedy search algorithm.



The proposed optimization framework is validated in an air quality sensing case study, based on real-world road geometry, sensing priorities, and cost structure. The results reveal considerable synergistic effect from fleet coordination, leading to not only greater sensing utilities, but also smaller gaps to the ideal sensing distributions. The mixed fleet sensing is also more cost-effective, with at least 30\% savings over any homogenous fleet. Furthermore, a transferability study based on three cities and hundreds of network variations highlight several key factors affecting the efficacy of mixed-fleet sensing, which serves as a practical guide to the planning of drive-by sensing infrastructure in urban areas.

The rest of this article is organized as follows. Section \ref{secLR} reviews related literature of drive-by sensing. The drive-by sensing coverage problem is formalized in Section \ref{secDSCprob}, followed by optimization models for the taxi+bus fleet (Section \ref{secTB}), DV fleet (Section \ref{secDV}), and taxi+bus+DV fleet (Section \ref{secTBD}). Section \ref{secCS}, presents a series of numerical experiments, and Section \ref{secconclude} offers some practical insights and future search directions.

\section{Related work}\label{secLR}

The deployment of mobile sensors is a nontrivial extension of the fixed sensor location problem \citep{VRCT2014, XHCC2016, ZFM2018, RCF2019, SKLLE2019}, due to the complex constraints introduced by the mobility patterns of the vehicles. The reader is referred to \cite{JHL2023} for a comprehensive review on the optimization aspect of drive-by sensing. Literature related to this study is presented in four parts: quantifying the sensing efficacy, and taxi-based, bus-based and DV-based drive-by sensing.

\subsection{Quantifying sensing efficacy}
Quantification of the sensing efficacy is foremost to the design and planning of drive-by sensing (DS) infrastructure. Existing literature on DS typically consider spatial objects such as road segments or spatial grids as the sensing target, and the sensing utility is quantified by:
\begin{enumerate}
\item a monotonically increasing utility function of the sensing frequency or data points in a time epoch \citep{PDB2012, ZMLL2015, ADFS2021}, or

\item a binary (0-1) utility, where $1$ is attained if and only if the subject is covered by at least $K$ times ($K\geq 1$) in a time epoch \citep{YDLCG2017, CDMFD2018, CCCFD2020, ZML2014}. 
\end{enumerate}

The former allows the utility to be accumulated through multiple visits, which is appropriate for objectives that require high sensing frequency (e.g. air quality, traffic conditions); The latter is most often employed with $K=1$, which is for objects with low temporal variations (e.g. road roughness, noise). In \cite{XCPJZN2019}, a prescribed sensing distribution is assumed and a given solution is evaluated by its distance from such a target distribution.

\subsection{Taxi-based drive-by sensing}
Taxis are popular sensor hosts in both research and practice, for their high spatial mobility and low costs. \cite{OASSR2019} employed a ball-in-bin model to quantify the sensing power of taxi fleet, and found that the sensing power of taxi fleet is large yet limited, as certain parts of the network are rarely visited. On the optimization taxi sensing power, one line of research optimizes sensor allocation through taxi subset selection with heuristic algorithms \citep{KGB2016, CLGZX2017, YFLRL2021, ZMLL2015, ZYL2018, HCL2015} or approximate algorithms \citep{YDLCG2017, LNL2016}. In another line, studies aim to exert influence on the operations of the fleet through incentivized routing for occupied \citep{ADFS2021, FJLQG2021} and unoccupied \citep{C2020, XCPJZN2019} taxis.

\subsection{Bus-based drive-by sensing}

Public transit vehicles such as buses and trams are also considered in DS for their extensive spatial coverage and regular operational patterns. \cite{CDMFD2018} use bus GPS data in Rio de Janeiro, and found a highly uneven distribution of coverage both spatially and temporally. Studies that optimize sensor deployment on bus fleets can be categorized according to their decision granularity, namely line-based allocation  \citep{Saukh2012, Kaivonen2020}, vehicle-based allocation  \citep{AD2017, TRMSO2020, WLQWL2018}, and trip-based allocation \citep{DH2023, JHL2023partC}.

\subsection{Drive-by sensing with dedicated vehicles}

Dedicated vehicles, with very high controllability and flexibility, are used to either perform highly targeted sensing \citep{JHG2023, LMZLL2023}, or to complement the sensing capabilities of other fleets \citep{FZGJGW2021}. The DSC problem with a fleet of DVs is related to the team orienteering problem \citep{DGM2013, KZLC2015, GM2020}, sometimes referred to as the prize-collecting vehicle routing problem \citep{LSPHZL2019, RHS2021}, in which the vehicles are coordinated to maximize the prizes collected at individual nodes or arcs within certain time or distance budgets. Nevertheless, the DSC problem for DVs considered in this work has two notable distinctions: (1) the `prizes' (i.e. sensing utilities) are collected at spatial grids instead of nodes/arcs while the vehicle trajectories must be fined to road networks, leading to a dual-spatial-scale routing structure; (2) the `prizes' can be accumulated by multiple vehicles because the sensing utility grows with more data collected. Such a new problem is termed dual-spatial-scale team orienteering problem with accumulable utilities.

\section{Problem statement and preliminary analyses}\label{secDSCprob}

\begin{definition}{\bf (The \underline{d}rive-by \underline{s}ensing \underline{c}overage (DSC) problem)}\label{defDSC}
Given a set of spatial units (grids or roads) as the subject of drive-by sensing, the DSC problem aims to allocate sensors to a fleet of vehicles of various types, and provide additional operational guidance (e.g. vehicle routing) as necessary, such that the spatial and temporal coverage is optimized. Such problems are often constrained with a fixed budget and diversified with customized sensing priorities and preferences. 
\end{definition}

Similar to facility location problems \citep{WWWZQ2021, RCS2020}, the DSC is strategic in nature since the sensors cannot be removed from the host vehicles. However, the fluidity of the sensors introduces temporal dimension to the problem, and the probabilistic nature of vehicle operations (such as taxis and buses) requires a mean-value approach. These features indicate that the DSC problem is a highly non-trivial extension of fixed sensor location problems found in the literature \citep{VRCT2014, XHCC2016, ZFZM2022}.

\subsection{Quantifying temporal-spatial coverage as sensing utilities}
The quantitative assessment of sensing utility is facilitated by appropriate space-time discretization. In this paper, space is meshed into grids $g\in \mathcal{G}$ while time is discretized as $t\in\mathcal{T}$. Let $N_{g,t}$ be the number of distinct vehicles whose trajectories intersect grid $g$ during time $t$, and $\xi\big(N_{g,t}\big)$ be the corresponding sensing utility. The following should hold for $\xi(\cdot)$:
\begin{itemize}
\item $\xi(\cdot)$, as a function of $N_{g,t}$, should be non-decreasing. In addition, as $N_{g,t}$ increases, the marginal utility decreases (i.e. the utility gain for one more vehicle coverage decreases). 
\end{itemize}

\noindent From a sampling point of view, $\xi(\cdot)$ should be non-decreasing. On the other hand, the diminishing marginal utility guarantees that no grid receives excessively high coverage, because it would be more efficient to allocate vehicle coverage to grids with lower utilities. Given such considerations, one viable choice for $\xi(\cdot)$ is: 
\begin{equation}\label{eqn1}
\xi\big(N_{g,t}\big)=\displaystyle\sum_{n=1}^{N_{g,t}} {1\over n^{\alpha}}\qquad \alpha\in(0,1),\quad N_{g,t}\in \mathbb{Z}_+
\end{equation}

\noindent It is desirable to extend the domain of $\xi(\cdot)$ to the set of non-negative real numbers $\mathbb{R}_+$, since we treat $N_{g,t}$ as a space-time-averaged quantity. To this end, we invoke the uniform discretization $0=n_1<n_2<\ldots<n_m=N_{g,t}$ with step size $\delta$, and write
\begin{equation}\label{rformal}
\xi(N_{g,t})\approx\sum_{m=1}^{N_{g,t}/\delta} {\delta\over n_m^{\alpha}}\approx \int_0^{N_{g,t}} {1\over x^{\alpha}}dx={1\over 1-\alpha}\big(N_{g,t}\big)^{1-\alpha}, \qquad N_{g,t}\in\mathbb{R}_+
\end{equation}
This leads to the following formula by introducing $\beta\doteq 1-\alpha$ and ignoring the multiplicative constant $1/(1-\alpha)$:
\begin{equation}\label{rformula}
\xi\big(N_{g,t}\big)=\big(N_{g,t}\big)^{\beta}\qquad\beta\in(0,1),\quad N_{g,t}\in\mathbb{R}_+
\end{equation} 
\noindent The utility function \eqref{rformula} will be used in the remainder of this paper. The choice of $\beta$ in various sensing context is discussed in Remark \ref{rmkbeta}. 

Building on the notion of sensing utility $\xi(N_{g,t})$, we further introduce spatial weights $w_g,\,g\in\mathcal{G}$ and temporal weights $\mu_t,\,t\in\mathcal{T}$, such that $\sum_{g\in\mathcal{G}}w_g=1$ and $\sum_{t\in\mathcal{T}}\mu_t=1$, to indicate the relative importance of information collected at $(g,t)$. These weights can be customized to accommodate various urban sensing needs.

\subsection{The DSC problem without vehicle mobility constraints}\label{subsecQSTC}

It is crucial to formulate and understand the objective function in the DSC problem to yield mathematically and physically meaningful solutions. This section provides an informative discussion of the optimization objective based on a conceptual DSC problem, by replacing complex vehicle mobility patterns with a simple resource constraint. Specifically, we consider the following problem:   
\begin{equation}\label{fdsceqn1}
\max_{{\bf N}=(N_{g,t}:\, g\in\mathcal{G}, t\in\mathcal{T})}\Phi({\bf N})=\sum_{t\in\mathcal{T}} \mu_t\sum_{g\in\mathcal{G}}w_g \xi\big(N_{g,t}\big)=\sum_{t\in\mathcal{T}} \sum_{g\in\mathcal{G}}\mu_t w_g \big(N_{g,t}\big)^{\beta}
\end{equation}
\begin{equation}\label{fdsceqn2}
\hbox{s.t.}~~h({\bf N})=\sum_{t\in\mathcal{T}}\sum_{g\in\mathcal{G}}N_{g,t}-M=0
\end{equation}
\noindent where $\Phi$ is the {\it \underline{s}pace-\underline{t}ime \underline{w}eighted \underline{s}ensing \underline{u}tility} (STWSU), which is the main objective function used throughout this paper. Eqn \eqref{fdsceqn2} is a conceptual resource constraint, which simply says that the total sum of vehicle coverage is bounded by $M$. This is a relaxed version of the DSC problem because we allow vehicle coverage to be entirely and freely splittable in space and time as if it were non-atomic fluid, while in reality it is confined to vehicle trajectories that follow specific mobility patterns (e.g. those of taxis, buses and DVs). Analyzing \eqref{fdsceqn1}-\eqref{fdsceqn2} could help us understand the target distribution of sensing utility in an ideal situation, and use it to benchmark those solutions with actual vehicle mobility constraints.

\subsubsection{First-order optimality conditions}\label{subsubsecKKT}
The KKT conditions of the convex program \eqref{fdsceqn1}-\eqref{fdsceqn2}, which is also sufficient, read
$$
\nabla\Phi({\bf N}^*) + \lambda \nabla h({\bf N}^*)=0~\Longrightarrow~\mu_tw_g\beta \big( N^*_{g,t}\big)^{\beta-1} =-\lambda \qquad \forall g\in\mathcal{G},\,t\in\mathcal{T}
$$ 
for some $\lambda$. This leads to the following:
\begin{equation*}
\xi\big(N^*_{g,t}\big)=\big(N^*_{g,t}\big)^{\beta}= \left(-\beta/ \lambda \right)^{\beta\over 1-\beta}(\mu_t w_g)^{\beta\over 1-\beta}\qquad \forall g\in\mathcal{G},\,t\in\mathcal{T}
\end{equation*}
In other words, in the optimal solution, we have that 
\begin{equation}\label{proptoeqn}
\xi\big(N^*_{g,t}\big) \propto (\mu_t w_g)^{\beta\over 1-\beta},\qquad N_{g,t}^*\propto \big(\mu_t w_g\big)^{1\over 1-\beta}\qquad \forall g\in\mathcal{G},\,t\in\mathcal{T}
\end{equation}
This leads to the following observations: 
\begin{enumerate}
\item By maximizing $\Phi$, the vehicle coverage tends to spread out spatially and temporally, as any pair $(g,\,t)$ with positive weights receives a positive utility, which is a desirable outcome in urban crowdsensing. 

\item In the optimal solution, pairs $(g,\,t)$ with higher weights $w_g$ and $\mu_t$ receive more coverage, which is intended by prescribing the spatial and temporal weights. 

\item As $\beta\to 0$, $\xi\big(N_{g,t}^*\big)$ tends to be uniformly distributed. This corresponds to cases where the marginal utility of repeated coverage is extremely low, and spatial-temporal extent outweighs local sensing frequency. Examples include urban heat island and road surface health. On the other hand, relatively large $\beta$ (between 0.2 and 0.5) encourages repeated sampling of important areas; examples include air quality and traffic congestion.

\item As $\beta\to 1$, the distributions of $\xi(N_{g,t}^*)$ and $N_{g,t}^*$ are singular, meaning the entire mass is concentrated in $(g,\,t)$ with the highest $\mu_t w_g$. This reinforces the idea that diminishing marginal utility ($\beta<1$) is needed to avoid over-concentration of sensing resources.  
\end{enumerate}

The parameter $\beta$ is crucial as it defines the utility function, which exert heavy influence on the drive-by sensing solution. While Observations 3-4 above offer intuitive insights in two extreme cases, it is not immediately clear how a proper $\beta$ should be chosen to meet specific sensing requirements. In remark \ref{rmkbeta} we offer a practical method to determine $\beta$.

\begin{remark}\label{rmkbeta}
 For simplicity, we ignore the temporal weights by setting: $\mu_t =1/|\mathcal{T}|,\forall t\in\mathcal{T}$, and obtain $N_{g}^*\propto (w_g)^{1/1-\beta}$. Let $g_1$ and $g_2$ be the grids with largest and smallest weights, respectively:
$$
w_{g_1}=\max_{g\in\mathcal{G}}w_g,\qquad w_{g_2}=\min_{g\in\mathcal{G}}w_g
$$
Then, one should find $\zeta>w_{g_1}/w_{g_2}$\footnote{This condition guarantees that $\beta\in(0, 1)$ per Eqn \eqref{betazeta}. In case $\zeta\leq w_{g_1}/w_{g_2}$, one should use a new set of weights $v_g\doteq (w_g)^a,\forall g\in\mathcal{G}$ where $a\in(0,1)$ is such that $\zeta>(w_{g_1}/ w_{g_2})^a$.}, such that $N_{g_1}^*=\zeta N_{g_2}^*$. In prose, one should specify the ratio $\zeta$ of vehicle coverage between the most and least important grids; such a ratio comes from domain knowledge and field practice.  Then, $\beta$ can be determined by invoking $N_{g}^*\propto (w_g)^{1/1-\beta}$:
\begin{equation}\label{betazeta}
\zeta={N_{g_1}^*\over N_{g_2}^*}=\left({w_{g_1}\over w_{g_2}}\right)^{1\over 1-\beta}\quad \Longrightarrow\quad \beta=1-\log_{\zeta}\left({w_{g_1}\over w_{g_2}}\right)
\end{equation}
For example, in the air quality sensing case study presented in Section \ref{secCS}, we set $\zeta=3$ (that is, vehicle number in the most important grid should triple that of the least important grid). In addition, $w_{g_1}/w_{g_2}=2.35$ (see Figure \ref{fignetwork}), which leads to $\beta=0.2$.  

Finally, \eqref{betazeta} implies that $\zeta$ increases with $\beta$, and therefore larger $\beta$ corresponds to more frequent sampling of important areas, which confirms Observation 3.
\end{remark}

\subsubsection{Target distribution}\label{subsubsecTD}

The optimality conditions give rise to the {\it \underline{t}arget \underline{d}istribution} of the utilities $\{\hbox{TD}_{g,t};\,g\in\mathcal{G}, t\in\mathcal{T}\}$, which can be compared with the {\it actual distribution} $\{\hbox{AD}_{g,t};\,g\in\mathcal{G}, t\in\mathcal{T}\}$ from a particular solution:
\begin{equation}\label{ADTD}
\hbox{TD}_{g,t}={\big(\mu_t w_g\big)^{\beta\over 1-\beta}\over \sum_{g', t'}\big(\mu_{t'} w_{g'}\big)^{\beta\over 1-\beta}}, \qquad \hbox{AD}_{g,t}={\xi(N_{g,t})\over \sum_{g',t'}\xi(N_{g',t'})},  \qquad\forall g\in\mathcal{G},~t\in\mathcal{T}
\end{equation}
\noindent Specifically, we use the Kullback-Leibler (KL) divergence $\hbox{KL}\big(\hbox{AD}||\hbox{TD}\big)$ from Bayesian inference \citep{KL1951}, which measures the change of information when one revises beliefs from the prior distribution $\hbox{TD}$ to the posterior distribution $\hbox{AD}$:
\begin{equation}\label{nsmg}
\hbox{KL}\big(\hbox{AD}|| \hbox{TD}\big)=\sum_{g\in \mathcal{G},t\in\mathcal{T}}\hbox{AD}_{g,t}\log{\hbox{AD}_{g,t}\over \hbox{TD}_{g,t}}
\end{equation}
Note that the target distribution is independent of the total resource $M$ in \eqref{fdsceqn2}. Given any DSC solution, the KL-divergence \eqref{nsmg} measures its deviation from an ideal distribution of sensing utilities. It can also indicate the loss of sensing efficacy caused by real-world vehicle mobility constraints.

In the remainder of this paper, we employ the STWSU $\Phi$ as the optimization objective, while using the KL-divergence as a second performance measure on an ex post facto basis.

\begin{remark}\label{rmkPTD}
Given a {\it prescribed target distribution} (PTD) $\{\hbox{PTD}_{g,t},\,g\in\mathcal{G}, t\in\mathcal{T}\}$ of vehicle coverage $N_{g,t}$, one may need to minimize the gap to such a PTD. In this case, we consider the following optimization problem: 
\begin{equation}\label{ficDSCphi}
\max_{{\bf N}=(N_{g,t}:\, g\in\mathcal{G}, t\in\mathcal{T})}\Phi({\bf N})=\sum_{t\in\mathcal{T}} \sum_{g\in\mathcal{G}}\mu_t w_g \big(N_{g,t}\big)^{\beta}=\sum_{t\in\mathcal{T}} \sum_{g\in\mathcal{G}}\big(\hbox{PTD}_{g,t}\big)^{1-\beta}\cdot \big(N_{g,t}\big)^{\beta}
\end{equation}
\begin{equation}\label{ficDSCcons}
\hbox{s.t.}~~h({\bf N})=\sum_{t\in\mathcal{T}}\sum_{g\in\mathcal{G}}N_{g,t}=M
\end{equation}
\noindent where we set $\mu_t w_g= \big(\hbox{PTD}_{g,t}\big)^{1-\beta}~ \forall g, t$. According to \eqref{proptoeqn}, the solution of this optimization problem satisfies:
$$
N^*_{g,t} \propto \big(\mu_t w_g\big)^{1\over 1-\beta}=\hbox{PTD}_{g,t}\qquad\forall g\in\mathcal{G},\,t\in\mathcal{T}
$$
\noindent This means that the distribution gap minimization problem \citep{XCPJZN2019} can be converted to a STWSU maximization problem by properly defining $\mu_t w_g$. 
\end{remark}


\subsection{The DSC problem with vehicle mobility constraints}

\begin{figure}[H]
\centering
\includegraphics[width=\textwidth]{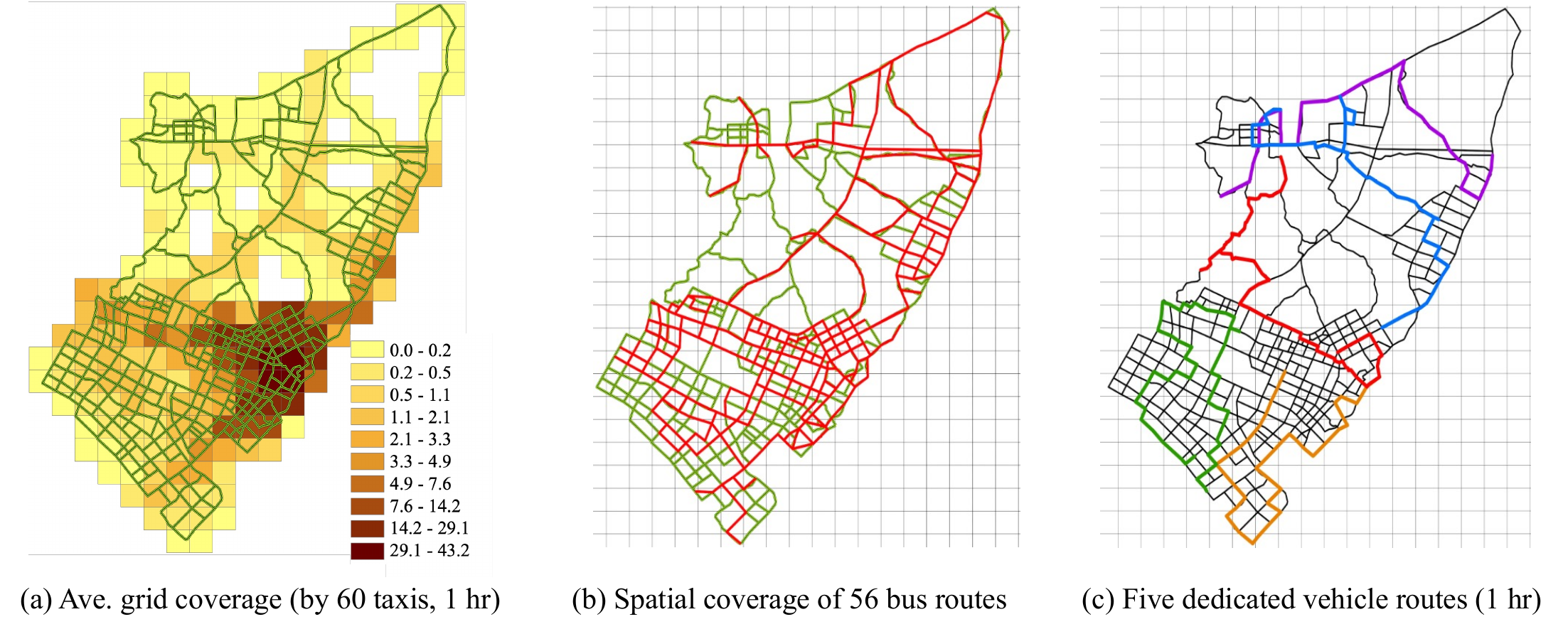}
\caption{Examples of the spatial mobility patterns of different fleet.}
\label{figTaxiBus}
\end{figure}

The optimization challenge of the DSC problem with mixed fleet lies in the heterogenous spatial-temporal mobility patterns of the host vehicles. In our setting, taxi-based sensing has low costs and long operating hours, but the spatial coverage is biased, primarily influenced by trip demand distribution; see Figure \ref{figTaxiBus}(a) for a real-world example. Transit vehicles such as buses have fixed routes, extending to remote areas not sufficiently covered by taxis, but could have substantial sensing blindspot (Figure \ref{figTaxiBus}(b)). Dedicated vehicles can be freely navigated in the network (Figure \ref{figTaxiBus}c), but the start-up, maintenance and personnel costs are considerably higher. 

The following assumptions are made on each vehicle type. 

{\bf Taxi} movements are stochastic, characterized by a ball-in-bin model \cite{OASSR2019}. Specifically, all the taxis are identically and independently distributed, and the number of taxis traversing a grid $g$ during time interval $t$ follows a binomial distribution $B(n^T, p_{g,t})$ where $n^T$ is the total number of instrumented taxis, and $p_{g,t}$ is the probability that a taxi covers grid $g$ during $t$. Under such a probabilistic model, the only decision variable of the taxi subproblem is $n^T$. 

{\bf Buses} are operated according to scheduled timetables. Each bus line $j$ has a fixed route and a number of trips to be served throughout the day. Given the bus GPS trajectories, one could estimate the en-route service time and turn-around time at terminals, from which the averaged number of buses covering a grid in one time interval can be estimated. The decision variable of the bus subproblem is the vector ${\bf y}^B=\big(y_1^B,\,\ldots,\,y_m^B\big)\in\mathbb{Z}_+^m$, where $y_j^B$ is the number of sensors allocated to line $j$, and $m$ is the total number of bus lines.

 {\bf Dedicated Vehicles:} The DVs can be freely navigated in the network. Each DV $i$ is assigned a fixed route $R_i^D$, to be repeatedly executed with round trips (other forms of DV operations are discussed in Section \ref{secmore}). The decision variables are ${\bf R}^D=\big\{R^D_1,\,\ldots,\,R^D_{n^D}\big\}$ where $n^D$ is the DV fleet size.

We may now present the following formulation of the DSC problem, where the mobility characteristics of different fleet types are presented as constraints.
\begin{align}
\label{DSCopt1}
&  \max_{n^T,\,{\bf y}^B,\,{\bf R}^D} \Phi = \sum_{t\in\mathcal{T}} \mu_t\sum_{g\in\mathcal{G}} w_g \big(N_{g,t}\big)^{\beta}
\\
\label{DSCopt2}
& N_{g,t}=N_{g,t}^T+N_{g,t}^B+N_{g,t}^D\quad \forall g,\,t
\\
\label{DSCopt3}
& N_{g,t}^T={\bf M}_{g,t}^T(n^T)\quad \forall g,\,t
\\
\label{DSCopt4}
& N_{g,t}^B={\bf M}_{g,t}^B({\bf y}^B)\quad \forall g,\,t
\\
\label{DSCopt5}
& N_{g,t}^D={\bf M}_{g,t}^D({\bf R}^D)\quad \forall g,\,t
\\
\label{DSCopt6}
&  c^T n^T +c^B \sum_{i=1}^my_i^B + c^Dn^D \leq M
\\
\label{DSCopt7}
&  n^T,\,y_j^B \in \mathbb{Z}_+,\quad n^T\leq L^T,\, y_j^B\leq L_j^B,\,j=1,\ldots,m
\end{align}

\noindent The objective \eqref{DSCopt1} is to maximize the STWSU. $N_{g,t}^T$, $N_{g,t}^B$, $N_{g,t}^D$ are, respectively, the expected number of taxis, buses and DVs traversing grid $g$ during $t$. \eqref{DSCopt3}-\eqref{DSCopt5} express the expected grid coverage by different vehicle types as functions/operators of their respective decision variables. \eqref{DSCopt6} is the budget constraint, where $c^T$, $c^B$, $c^D$ are the unit costs for instrumenting a taxi, bus and dedicated vehicle, which includes sensor procurement/maintenance cost as well as installation and vehicle operation costs (a detailed instance of the cost structure is presented in Table \ref{tabcs}). \eqref{DSCopt7} imposes the total taxi fleet size $L^T$ and total number of buses per line $L_j^B$ as upper bounds.

\section{The taxi-bus subproblem}\label{secTB}
This section aims to solve the DSC problem using taxis and buses as host vehicles. 

\subsection{Approximation of taxi coverage}\label{subsecAppT}
As suggested by \cite{OASSR2019}, the coverage of a grid by taxis can be characterized by a binomial distribution. Indeed, assuming $n^T$ taxis are operating independently, following identical spatial distribution in the network, the number of taxis that cover a grid $g$ within $t$ is binomially distributed $m_{g,t} \sim B(n^T,\,p_{g,t})$, where $p_{g,t}$ is the probability that a single taxi intersects $g$ during time $t$. Therefore, the average number of taxis covering $g$ during $t$ can be expressed as the expectation of the binomial distribution: 
\begin{equation}\label{MT}
N_{g,t}^T={\bf M}_{g,t}^T(n^T) \approx n^T p_{g,t}\qquad g\in \mathcal{G},\,t\in\mathcal{T}
\end{equation}

\begin{figure}[h!]
\centering
\includegraphics[width=.75\textwidth]{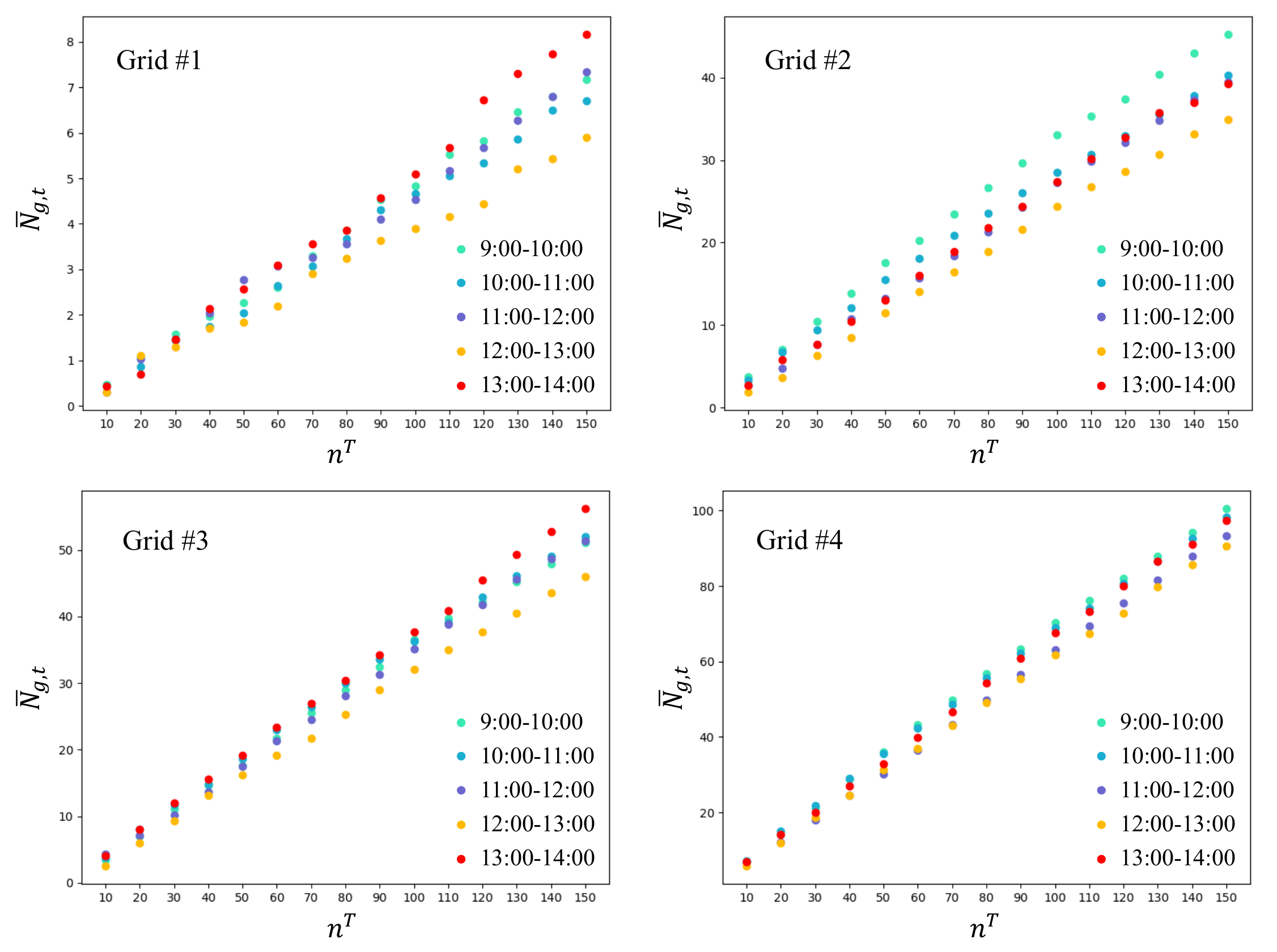}
\caption{Empirical relationship between fleet size $n^T$ and average coverage $\bar N_{g,t}$.}
\label{figBinomial}
\end{figure}

Next, we use empirical data from Longquanyi District (Chengdu, China) to validate the approximation \eqref{MT}. We consider 5 hrs (9am-14pm) in 31 days of March 2021. For a fixed fleet size $n^T$, let $N_{g,t}(\tau)$ be the number of taxis covering grid $g$ during hr $t$ ($t=9,\ldots,13$) on the $\tau$-th day ($\tau=1,\ldots,31$), and define 
\begin{equation}\label{Nbar}
\bar N_{g,t}={1\over 31}\sum_{\tau=1}^{31}N_{g,t}(\tau)\qquad g\in \mathcal{G},\quad t=9,\ldots, 13
\end{equation}
\noindent to be the monthly average, which reduces daily variations, given the planning nature of the DSC problem. The empirical relationships between $n^T$ and $\bar N_{g,t}$ are illustrated in Figure \ref{figBinomial}. We see that, for each hour, the relationships are approximately linear, which agrees with \eqref{MT}. In addition, the slope $p_{g,t}$ of such a linear form varies not only by the hour but also across different grids.


\begin{figure}[H]
\centering
\includegraphics[width=\textwidth]{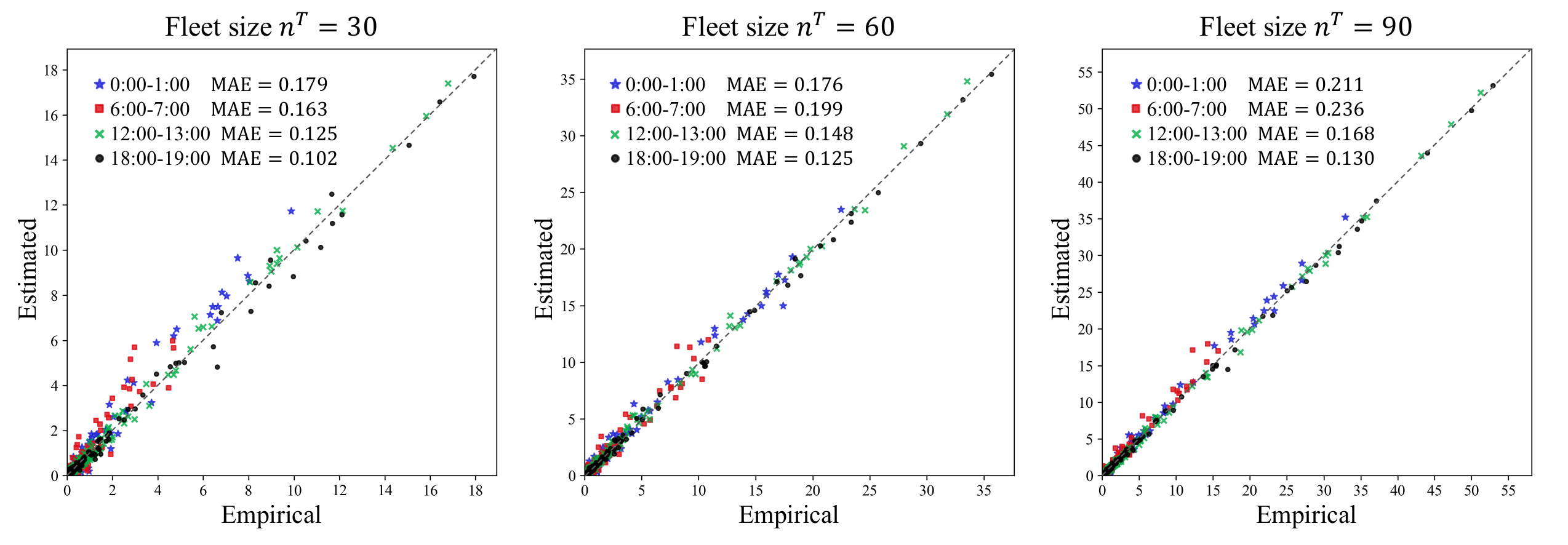}
\caption{Scatter plots of empirical vs. estimated coverage for all $g\in\mathcal{G}$ (221 grids), during four different hours. Each sub-figure corresponds to a given fleet size $n^T$.}
\label{figbinomialt}
\end{figure}

Then, the parameters $p^*_{g,t}$'s are obtained by fitting the linear relationship $\bar N_{g,t}\sim n^Tp_{g,t}$, where multiple random draws of $n^T$ taxis are performed to estimate the mean value of $\bar N_{g,t}$, $\forall g\in\mathcal{G},\,1\leq t\leq 24$. The purpose of the random draws is to reduce the selection bias. In Figure \ref{figbinomialt} we compare the empirical grid coverage $\bar{N}_{g,t}$ with the ones estimated as $p^*_{g,t}n^T$, for $g\in\mathcal{G}$ and $t=1, 7, 13, 19$ (other times are omitted due to space limitation). Firstly, we observe that all the scatter points are centered around the $y=x$ axis, suggesting the validity of the binomial approximation. Secondly, for the same $n^T$, the mean absolute errors (MAEs) are smaller for 12:00-13:00 \& 18:00-19:00 due to the law of large numbers, as taxis are most active during these times (see Figure \ref{figflow}). In particular, the errors are very low for high-volume grids.

\subsection{Approximation of bus coverage}\label{subsecAppB}

The duration of a one-way bus trip on line $j$ consists of an en-route service time $T^s_j$ (hrs) and a turn-around time $T^a_j$ (hrs) at bus terminals (for headway control, fueling/charging, staff change, etc.). These times also change over the course of a day due to dynamic traffic congestion and time-varying service frequencies. The quantities $T_j^s(t)$ and $T_j^a(t)$ can be estimated using bus operational data (e.g. GPS trajectories or timetables). 

\begin{figure}[h!]
\centering
\includegraphics[width=.7\textwidth]{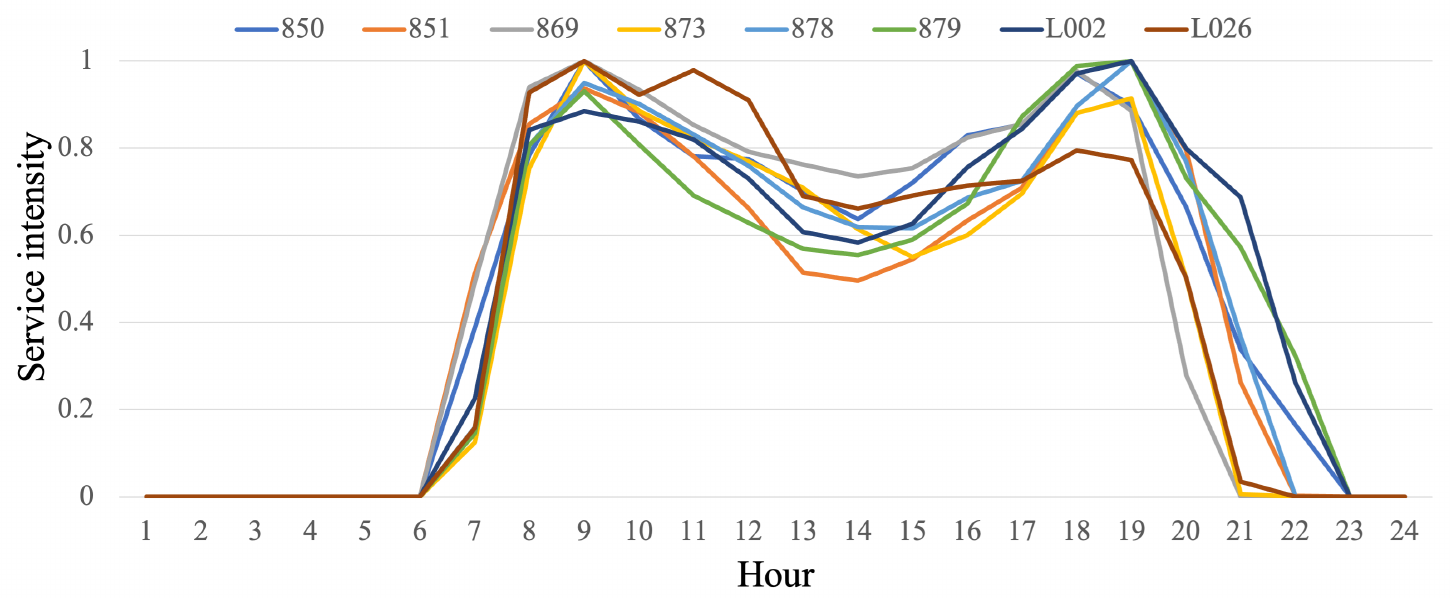}
\caption{Service intensity $\gamma_j(t)$ of a few bus lines in Longquanyi District.}
\label{figserint}
\end{figure}

Each bus line $j$ has its own fleet of size $L_j^B$, among which $y_j^B$ buses are instrumented ($y_j^B$ is the decision variable here). Note that not all the $L_j^B$ buses are in service during any time of the day. In fact, let $\lambda_j(t)$ be the average number of buses in service during time $t$, and define the service intensity:
\begin{equation}\label{gammaeqn}
\gamma_j(t)={\lambda_j(t)\over  L_j^B}\in [0,\,1], \qquad t\in\mathcal{T}
\end{equation}
\noindent $\gamma_j(t)$ can be estimated from timetables or bus GPS data. The service intensities of a few selected lines in Longquanyi District are shown in Figure \ref{figserint}. Under the premise that bus operation is independent from sensor instrumentation, each bus has an equal probability $\gamma_j(t)$ of being in service. Therefore, the expected number of instrumented buses in service during time $t$ is $\gamma_j(t)y_j^B$. 

Again, given the strategic nature of the problem, we employ a time-aggregated approach. At time $t$ the period for serving a one-way trip is $T_j^s(t)+T_j^a(t)$, thus each grid on the bus route is traversed $1/(T_j^s(t)+T_j^a(t))$ times per hr by one bus. The number of sensing coverage is expressed as $\gamma_j(t)y_j^B/(T_j^s(t)+T_j^a(t))$.

We now consider an arbitrary grid $g\in \mathcal{G}$, and define
$$
\delta_{gj}=
\begin{cases}
1\quad &\text{if bus route}~j~\text{intersects grid}~g
\\
0\quad &\text{otherwise} 
\end{cases},
$$
\noindent The total number of coverage during time $t$ is consequently expressed as 
\begin{equation}\label{MB}
N_{g,t}^B={\bf M}_{g,t}^B({\bf y}^B) =  \sum_{j=1}^m\delta_{gj}{\gamma_j(t)y_j^B\over T_j^s(t)+T_j^a(t)}\qquad g\in \mathcal{G},\,t\in\mathcal{T}
\end{equation}

\begin{remark}
The proposed bus subproblem is actually a probabilistic and time-averaged approach, neglecting details such as dispatch schedules, vehicle movement, and operational variations. In another perspective, we treat the bus fleet as a continuum, circumventing the attempt to model individual bus operations (which is difficult as it requires much more for a planning problem). This is consistent with our definition of sensing utility $\xi(N)$, where $N$ is treated as a real number.
\end{remark}

\subsection{Convex formulation of DSC with taxis and buses}\label{subsecDSCTB}

With the mobility constraints for taxis and buses instantiated as  \eqref{MT} and \eqref{MB}, the DSC problem with taxi-bus fleet is expressed as follow. 
\begin{equation}\label{pconvext1}
\max_{n^T,\,y_1^B,\,\ldots,\,y_m^B}\Phi= \sum_{t\in\mathcal{T}} \mu_t \sum_{g\in\mathcal{G}} w_g \left(n^Tp_{g,t} + \sum_{j=1}^m \delta_{gj}{\gamma_j(t) y_j^B\over T_j^s(t)+T_j^a(t)} \right)^{\beta}
\end{equation}
\begin{equation}\label{pconvext2}
c^T n^T + c^B\sum_{j=1}^m y_j^B\leq M
\end{equation}
\begin{equation}\label{pconvext3}
0\leq n^T\leq L^T,~~ 0\leq y_j^B\leq L_j^B,~~j=1,\ldots,m
\end{equation}
\noindent where $\beta\in (0,\,1)$; $p_{g,t}$'s are estimated based on empirical taxi GPS data; $\delta_{gj}$'s are obtained from bus route geometry; $\gamma_j(t)$, $T_j^s(t)$ and $T_j^a(t)$ are obtained from bus GPS data; $L^T$ is the total number of taxis and $L_j^B$ is the total number of buses in line $j$. Note that the decision variables $n^T$ and $y_j^B$ are relaxed to be real numbers in \eqref{pconvext3}, so that the optimization problem becomes convex. Its real-valued solution can be easily converted to integers values that satisfy the constraints. 

As a side note, the formulation \eqref{pconvext1}-\eqref{pconvext3} can be used to optimize single-mode DSC problems by setting $L_j^B=0, \forall j$ (taxi only) or $n^T=0$ (bus only).

\begin{remark}
In the objective \eqref{pconvext1}, grid coverage contributed by vehicles other than taxis or buses can be modeled with an additive constant inside $(\cdot)^{\beta}$, without affecting the convexity w.r.t. $n^T$ and $y_j^B$'s. This is particularly important for the joint DSC problem with taxis, buses and dedicated vehicles to be presented later. 
\end{remark}

\section{The dedicated vehicle subproblem}\label{secDV}

Unlike third-party vehicles (taxi, bus), dedicated vehicles (DVs) can be freely navigated in the road network to perform data collection tasks. The corresponding DSC problem resembles the {\it team orienteering problem} \citep{DGM2013, KZLC2015}, with two notable distinctions: (1) the sensing `scores' are collected by intersecting grids instead of visiting nodes or arcs; (2) the scores can be accumulated according to the sensing utility function $\xi(\cdot)$. 

In light of the following discussion, the DSC problem with DVs is termed {\it dual-spatial-scale team orienteering problem with accumulable utilities}. The key challenge is to coordinate grid coverage with network routing, as they are two codependent processes taking place on distinct spatial scales (Definition \ref{defNGroute}). Another issue is the fact that the utility obtained by traversing a grid is not always a constant (Definition \ref{defMU}).

Let $\mathcal{N}$ be the set of nodes in the road network, and $\mathcal{G}$ be the set of grids. If a node $v\in\mathcal{N}$ is located in a grid $g\in\mathcal{G}$, we write $v\in g$. 


\begin{definition}\label{defNGroute}{\bf (N-route and G-route)}
A N-route (N stands for network) refers to a conventional route in the road network, expressed as a sequence of nodes $r_N=(v_1,\,\ldots,\,v_m)\subset \mathcal{N}$; a G-route (G stands for grid) refers to a conceptual route in the grid space, expressed as the sequence of grids $r_G=(g_1,\,\ldots,\,g_n)\subset\mathcal{G}$ being traversed. The sets of N-routes and G-routes are respectively denoted $\mathcal{R}_N$ and $\mathcal{R}_G$.
\end{definition}

\begin{definition}\label{defMU}{\bf (Marginal utility)}
For a grid $g\in\mathcal{G}$ with sensing weight $w_g$, which is covered by $N_g$ sensing vehicles, the marginal sensing utility is defined as $\eta_g=\eta_g(N_g)=w_g\big(\xi (N_g+1)-\xi(N_g)\big)$. 
\end{definition}

\begin{figure}[h!]
\centering
\includegraphics[width=\textwidth]{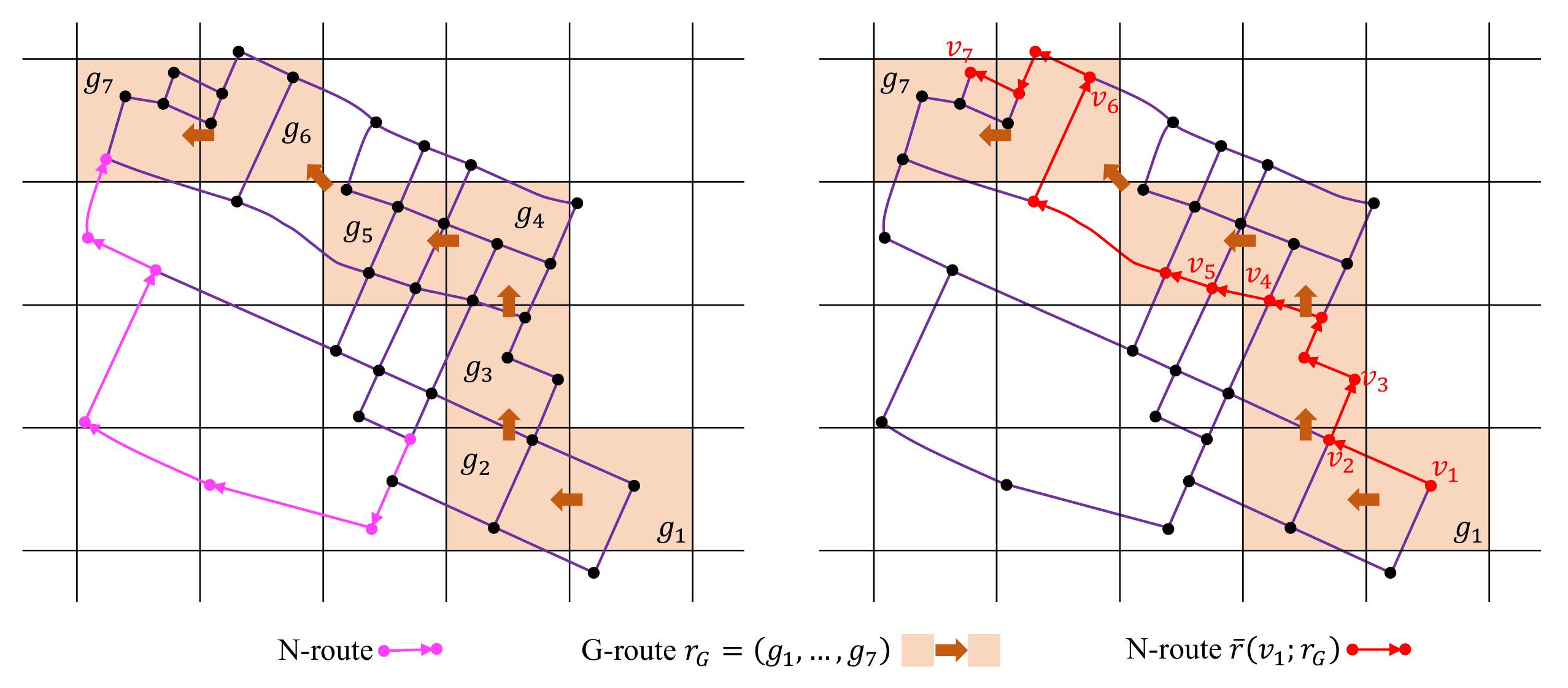}
\caption{Example of the dual-spatial-scale structure. The left figure shows a N-route and a G-route; the right figure illustrates the N-route that executes the G-route with starting node $v_1\in g_1$. Note that the G-route traversed by $\bar r(v_1; r_G)$ is different from $r_G$.}
\label{figNGroutes}
\end{figure}

It is natural to assume that any grid subject to mobile sensing can be reached from any given point in the road network. Next, we elaborate the relationship between a N-route and G-route. Given any G-route $r_G=(g_1,\,\ldots,\,g_n)$, and a starting node $v_1\in g_1$, the N-route that executes the movement $g_1 \to g_2$ needs to start with $v_1$ and ends with the closest node $v_2\in g_2$, where closeness is measured by network distance. Using $v_2$ as the starting point, the movement $g_2\to g_3$ is executed in the same way. Therefore, the N-route that fully executes $r_G=(g_1,\,\ldots,\,g_n)\in\mathcal{R}_G$ with starting node $v_1\in g_1$, denoted $\bar r(v_1;\,r_G)\in\mathcal{R}_N$, is obtained by concatenating these intermediate legs $v_1\to v_2$, $v_2 \to v_3$, \ldots $v_{n-1}\to v_{n}$. Such a procedure is formalized in Algorithm \ref{algroute} and illustrated in Figure \ref{figNGroutes}. Note that, by following N-route $\bar r(v_1; r_G)$, the G-route being actually traversed is denoted $G_{\bar r(v_1; r_G)}$, which may be different from $r_G$.

\begin{algorithm}[h!]
 \caption{(Executing a G-route)}
\begin{tabbing}
\hspace{0.01 in}\=  \hspace{0.9 in}\= \kill 
\>{\bf Input}  \> G-route $r_G=(g_1,\ldots,g_n)$, starting node $v_1\in g_1$.
\\
  \>  {\bf Step 1} \>  For $i=1,\ldots, n-1$
  \\
  \> \> \hspace{0.2 in} Find $v_{i+1}=\underset{v\in g_{i+1}}{\text{argmin}}\,dist(v,\,v_i)$, and the shortest path $p_i^*=(v_i,\ldots, v_{i+1})$
 \\
 \> \> End For
 \\
 \>  {\bf Step 2}   \> Concatenate $p_i^*$'s to form the full N-route $\bar r(v_1; r_G)=p_1^*\oplus p_2^*\ldots \oplus p_{n-1}^*$
 \\
 \>{\bf Step 3} \> Calculate the actual G-route $G_{\bar r(v_1; r_G)}$ traversed by $\bar r(v_1; r_G)$.
 \\
\>{\bf Output}      \> N-route $\bar r(v_1; r_G)$, destination node $v_n$, G-route $G_{\bar r(v_1; r_G)}$.
       \end{tabbing}
 \label{algroute}
\end{algorithm}

\subsection{Dual-spatial-scale routing for a single DV}

The routing of a single DV has a fixed time window $[T_1,\,T_2]$ and an origin node $o$ in the network. The overall strategy is to conduct a local search of destination grid, as well as a G-route to that grid, to collect sensing utilities in a time-efficient way. This process will be repeated, while updating the marginal utilities of traversed grids to avoid repeated visits, until the time budget $T_2-T_1$ is reached.


Given the origin node in the road network $o\in g_o\in \mathcal{G}$, the idea is to find a destination grid $g_d^*$ and a G-route that maximize the utility efficiency $\psi(\cdot)$:
\begin{equation}\label{psieqn}
g_d^*=\underset{g_d}{\text{argmax}}\,\psi(g_d)=\underset{g_d}{\text{argmax}}{\sum_{g\in G_{\bar r(o;\,r_G)}} w_g \eta_g\over \phi\big(\bar r(o;\,r_G)\big)}
\end{equation}

\noindent where $r_G$ is the least-impeded G-route from $g_o$ to $g_d$, which is obtained from a modified A-star algorithm where the impedance of a grid $g$ is the inverse of its marginal utility $1/\eta_g$. Such a configuration of the A-star algorithm encourages the DV to traverse grids with high marginal utilities. $\phi\big(\bar r(o;\,r_G)\big)$ denotes the total travel time along the N-route $\bar r(o;\,r_G)$. The utility efficiency $\psi(g_d)$ is the ratio between the weighted sensing utilities obtained along the G-route $G_{\bar r(o;\,r_G)}$, and the route travel time. By maximizing $\psi(g_d)$, the routing strategy aims to collect more utilities within shorter time.

Given that the objective function $\psi(\cdot)$ has embedded A-star search and dual-spatial-scale routing components, which is highly complex, we conduct a destination search by enumerating all grids within a certain radius $R$ (e.g. 5 km) from the origin grid $g_o$, as outlined in Algorithm \ref{alg2}.

\begin{algorithm}[H]
 \caption{(Generation of local destination and route)}
\begin{tabbing}
\hspace{0.01 in}\=  \hspace{0.8 in}\= \kill 
\>{\bf Input}  \> Origin node $o\in g_o\in \mathcal{G}$, marginal utilities $\{\eta_g,\,g\in\mathcal{G}\}$
\\
  \>  {\bf Step 1} \>  Find the set of all grids within radius $R$ from the origin, denoted $\mathcal{G}_F$;
 \\
 \>  {\bf Step 2}   \> For $g_d\in \mathcal{G}_F$
 \\
 \>    {\bf Step 2a}        \>  \hspace{0.2 in} With origin $g_o$ and destination $g_d$, apply the modified A-star algorithm using  
 \\
 \> \> \hspace{0.2 in} the inverse of the marginal utility as the grid impedance, to obtain $r_G$;
 \\
 \>{\bf Step 2b} \> \hspace{0.2 in} Find $\bar r(o;\,r_G)$ and $G_{\bar r(o;\,r_G)}$ using Algorithm \ref{algroute}; calculate  $\psi(g_d)$ using \eqref{psieqn};
 \\
 \>                     \> End For
 \\
  \> {\bf Step 3}   \> Find $g^*_d=\underset{g_d\in\mathcal{G}_F}{\text{argmax}}\psi(g_d)$.
\\
\>{\bf Output}      \> Destination $d^*\in g_d^*$, G-route $r_G^*$ from $g_o$ to $g_d^*$, N-route $\bar r(o;\,r_G^*)$.
       \end{tabbing}
 \label{alg2}
\end{algorithm}


The generation of a complete route within time window $[T_1,\,T_2]$ follows Algorithm \ref{alg3}, by iteratively executing Algorithm \ref{alg2} while updating the marginal utilities of relevant grids.

\begin{algorithm}[H]
 \caption{(Single DV routing with time budget)}
\begin{tabbing}
\hspace{0.01 in}\=  \hspace{0.8 in}\= \kill 
\>{\bf Input}  \> Time window $[T_1,\,T_2]$, starting node $s\in g_s$, marginal utilities $\{\eta_g,\,g\in\mathcal{G}\}$
\\
\> {\bf Initialize} \> $o=s$, $g_o=g_s$; total travel time $\mathcal{\phi}=0$.
\\
  \>  {\bf Step 1} \>  With $o\in g_o$, apply Algorithm \ref{alg2} to find destination $d^*\in g_d^*$, and $\bar r(o;\,r_G^*)$;
 \\
 \> \> Update marginal utilities $\{\eta_g,\, g\in G_{\bar r(o;r_G^*)}\}$;
 \\
 \>     \> Calculate route travel time $TT(\bar r(o;\,r_G^*))$, let $\phi=\phi+TT(\bar r(o;\,r_G^*))$;
 \\
  \> {\bf Step 2}   \> If $T_1+\phi<T_2$, let $o=d^*$, $g_o=g_d^*$, and go to Step 1; 
  \\
  \>                      \> otherwise, find the node within $\bar r(o;\,r_G^*)$ whose arrival time is the closest to $T_2$,
  \\
  \>                      \> terminate the algorithm;
\\
\>{\bf Output}      \> Actual N-route $R^D$ followed within $[T_1,\,T_2]$, total weighted sensing utilities, marginal 
\\
\> \> utilities $\{\eta_g,\, g\in\mathcal{G}\}$.
       \end{tabbing}
 \label{alg3}
\end{algorithm}

\subsection{Dual-spatial-scale routing for multiple DVs}

The routing for multiple DVs needs to consider potential overlap of their search spaces. We employ a sequential optimization procedure with adjustment as outlined in Algorithm \ref{alg4}.

\begin{algorithm}[H]
 \caption{(Multi-vehicle routing within $[T_1,\,T_2]$)}
\begin{tabbing}
\hspace{0.01 in}\=  \hspace{0.8 in}\= \kill 
\>{\bf Input}  \> Time window $[T_1,\,T_2]$, starting nodes of the vehicles $s_n,\,n=1\ldots,n^D$,  marginal 
\\
\> \> utilities $\{\eta_g,\,g\in \mathcal{G}\}$, iteration limit $\text{max}\_\text{iter}$.
\\
\> {\bf Initialize} \> Set $k=1$
\\
  \>  {\bf Step 1} \>  Solve single-vehicle routing problems sequentially:
  \\
  \>      \> For $n=1,\ldots,n^D$
  \\
  \>     \> \hspace{0.2 in} Apply Algorithm \ref{alg3} to the $n$-th vehicle with marginal utilities $\{\eta_g,\,g\in\mathcal{G}\}$;
 \\
 \>       \> End For
 \\
 \>  {\bf Step 2}   \> Adjust routing solutions: 
 \\
 \>      \> For $n=1,\ldots, n^D$
 \\
 \>  {\bf Step 2a}   \> \hspace{0.2 in} Calculate the marginal utilities $\{\eta_g^{(\bar n)},\,g\in\mathcal{G}\}$ by removing the coverage of 
 \\
 \> \> \hspace{0.2 in} the $n$-th vehicle;
 \\
 \> {\bf Step 2b}\> \hspace{0.2 in} Apply Algorithm \ref{alg3} to the $n$-th vehicle with $\{\eta_g^{(\bar n)},\,g\in\mathcal{G}\}$, if the new solution 
 \\
 \> \> \hspace{0.2 in} yields higher utilities than the current solution, then accept the new solution 
 \\
 \> \> \hspace{0.2 in} and update $\{\eta_g,\,g\in\mathcal{G}\}$; otherwise, continue with $n=n+1$;
 \\
 \> \> End For
 \\
  \> {\bf Step 3}   \> If no solutions are updated for any of the vehicles or $k=\text{max}\_\text{iter}$, terminate  
  \\
  \> \> the algorithm; otherwise, set $k=k+1$ and go to Step 2.
\\
\>{\bf Output}      \>  Routing solutions ${\bf R}^D=\big\{R_1^D,\ldots, R_{n^D}^D \big\}$ for $n^D$ DVs.
       \end{tabbing}
 \label{alg4}
\end{algorithm}

Recall that the sensing utilities $\xi(N_{g,t})$ are calculated for every discrete time window $t\in\mathcal{T}$. Therefore, it is sensible to plan the DV routes within such time window $t$, and reverse-travel these routes during $t+1$. In this way, we not only ensure the same set of grids are covered for each $t\in \mathcal{T}$, but also avoid the need to plan long DV routes spanning several time windows, which is computationally expensive. This is also favorable in practice because assigning fixed route/area for periodic patrol not only reinforces the drivers' familiarity with local environment, but also allows flexible route adjustment in the event of emergencies. 

In light of the above discussion, we extend the routing plan from Algorithm \ref{alg4} to the entire time horizon $\mathcal{T}$, by setting $[T_1,\,T_2]$ as one single time window $t$, and turning the routes $R_1^D,\ldots, R_{n^D}^D$ into repeatable round trips, such that a one-way trip is performed in every $t$.

Finally, unlike traditional vehicle routing problems with given depots as starting locations, the DSC problem allows the DV routes to start at any node in the network, because the trip to/from the depot is negligible as far as sensing is concerned. In addition, although the choice of the starting locations may impact the overall sensing performance, it is not the focus of this work. Figure \ref{figTaxiBus}(c) shows five DV routes within an hour, with different starting locations, which are generated by the proposed algorithm.

\section{The DSC problem with three vehicle types}\label{secTBD}

\subsection{Problem formulation}

The three vehicle types (taxi, bus, dedicated vehicle) are simultaneously considered in the DSC problem. This is formulated as below.
\begin{equation}\label{DSCfeqn1}
\max_{\substack{n^T \\ {\bf y}^B=(y_1^B,\ldots,y_m^B) \\ {\bf R}^{D}=\{R_1^D,\ldots, R_{n^D}^D\}}} \Phi =\sum_{t\in\mathcal{T}}\mu_t\sum_{g\in\mathcal{G}}w_g \Big(N_{g,t}^T + N_{g,t}^B + N_{g,t}^D \Big)^{\beta}
\end{equation}
\begin{equation}\label{DSCdettaxiMP}
N_{g,t}^T=n^Tp_{g,t}\qquad\forall g\in\mathcal{G},\,t\in\mathcal{T}
\end{equation}
\begin{equation}\label{DSCdetbusMP}
N_{g,t}^B=\sum_{j=1}^m \delta_{gj}{\gamma_j(t) y_j^B\over T_j^s(t)+T_j^a(t)}\quad\forall g,\,t
\end{equation}
\begin{equation}\label{DSCdetDVMP}
N_{g,t}^D=\sum_{k=1}^{n^D}\sigma_{g,t}^k \quad\forall g,\,t
\end{equation}
\begin{equation}
c^T n^T +c^B \sum_{i=1}^my_i^B + c^Dn^D \leq M
\end{equation}
\begin{equation}\label{DSCfeqn6}
0\leq n^T\leq L^T,\, 0\leq y_j^B\leq L_j^B,\,j=1,\ldots,m
\end{equation}
\noindent where
\begin{equation}\label{defsigma}
\sigma_{g,t}^k=\begin{cases}
1\quad \hbox{if}~R_k^D ~\hbox{traverses}~ g~ \hbox{during}~ $t$
\\
0\quad \hbox{otherwise}
\end{cases},\quad k=1,\ldots, n^D
\end{equation}

The formulation above instantiates the conceptual DSC problem \eqref{DSCopt1}-\eqref{DSCopt7} by explicitly expressing the mobility constraints for taxis \eqref{DSCdettaxiMP}, buses \eqref{DSCdetbusMP}, and DVs \eqref{DSCdetDVMP}.

\subsection{Solution algorithm}

Given that the taxi-bus subproblem is convex (exact) and the DV algorithm is heuristic, we adopt a decomposition scheme that iterates between these two subproblems. More specifically, let $M^{T+B}$ be the budget for the taxi-bus subproblem, and $n^D$ be the number of DVs in the DV subproblem. Clearly $M^{T+B}+c^Dn^D\leq M$. We write:
\begin{equation}\label{blpuln1}
\text{Taxi-bus subproblem:}~~(n^{T},\,{\bf y}^{B})=\Psi_1\big(M^{T+B};\, {\bf R}^D \big)
\end{equation}
\begin{equation}\label{blpuln2}
\text{DV subproblem:}~~{\bf R}^{D}=\Psi_2\big(n^D;\,n^T,\,{\bf y}^B\big)
\end{equation}
\noindent where 
$$
\Psi_1\big(M^{T+B};\,{\bf R}^D\big)
$$ 
denotes the solution of the convex program \eqref{pconvext1}-\eqref{pconvext3} with budget $M^{T+B}$ and a priori coverage afforded by DVs (${\bf R}^D$). Symmetrically, 
$$
\Psi_2\big(n^D;\,n^T,\,{\bf y}^B\big)
$$ 
denotes solution of the DV routing problem provided by Algorithm \ref{alg4} with $n^D$ vehicles, and marginal utilities $\{\eta_g,\,g\in\mathcal{G}\}$ based on the existing taxi-bus coverage ($n^T$, ${\bf y}^B$).

The subproblems \eqref{blpuln1} and \eqref{blpuln2} are coupled via the budget constraint and shared spatial-temporal coverage. Algorithm \ref{alg5} proposes a heuristic method to solve both problems simultaneously by generating the following intermediate iterates: ($k=1,2,\ldots$)
\begin{equation}\label{subpiter1}
(n^{T, k},\,{\bf y}^{B,k})=\Psi_1\big(M^{T+B};\, {\bf R}^{D,k} \big)
\end{equation}
\begin{equation}\label{subpiter2}
{\bf R}^{D,k+1}=\Psi_2\big(M^D;\,n^{T, k},\,{\bf y}^{B, k}\big)
\end{equation}

\begin{algorithm}[h!]
 \caption{(DSC problem with taxi-bus-DV fleet)}
\begin{tabbing}
\hspace{0.01 in}\=  \hspace{0.8 in}\= \kill 
\>{\bf Input}  \> Total budget $M$, cost coefficients $c^T$, $c^B$, $c^D$, iteration limit `$\text{max}\_\text{iter}$'
\\
\> \> For $n^D=0: \lfloor M/c^D\rfloor$
\\
\> {\bf Initialize} \> \hspace{0.2 in} Let $k=0$; objective $\Phi^k=0$; set the DV route set ${\bf R}^{D, k}=\emptyset$;
\\
  \>  {\bf Step 1} \> \hspace{0.2 in} With budget $M^{T+B}=M-c^D n^D$, solve the taxi-bus subproblem  
  \\
  \> \> \hspace{0.2 in}  $\Psi_1\big(M^{T+B};{\bf R}^{D,k} \big)$ to obtain $n^{T,k}$ and ${\bf y}^{B,k}$;
  \\
  \>  {\bf Step 2}    \> \hspace{0.2 in} Update routes ${\bf R}^{D,k}=\Psi_2\big(n^D; \, n^{T,k},{\bf y}^{B,k}\big)$ by applying Algorithm \ref{alg4};
  \\
  \>{\bf Step 3} \> \hspace{0.2 in} Calculate the objective $\Phi^{k+1}$ based on $n^{T,k}$, ${\bf y}^{B,k}$ and ${\bf R}^{D,k}$. If $\Phi^{k+1}\leq \Phi^k$  
  \\
  \> \> \hspace{0.2 in}  or $k=\text{max}\_\text{iter}$, let $\big(\Phi^*_{n^D},\,n^{T,*}_{n^D},\,{\bf y}^{B,*}_{n^D},\,{\bf R}^{D,*}_{n^D}\big)=\big(\Phi^{k},\, n^{T,k},\,{\bf y}^{B,k},\,{\bf R}^{D,k}\big)$ and 
  \\
  \> \> \hspace{0.2 in}  continue For; Otherwise, let ${\bf R}^{D,k+1}={\bf R}^{D,k}$, set $k=k+1$ and go to Step 1.
  \\
  \> \> End For
  \\
 \> {\bf Step 4}      \>  Let $n^{D,*}=\underset{0\leq n^D\leq \lfloor M/c^D\rfloor}{\text{argmax}}\Phi_{n^D}^*$
\\\\
\>{\bf Output}      \>  DSC solution $\big(n^{T,*}_{n^{D,*}} ,\,{\bf y}^{B,*}_{n^{D,*}},\,{\bf R}^{D,*}_{n^{D,*}}\big)$, with objective $\Phi^*_{n^{D,*}}$.
       \end{tabbing}
 \label{alg5}
\end{algorithm}

In Algorithm \ref{alg5}, the outer loop (`For') enumerates budget split between taxi-bus and DVs, while the inner iteration (Steps 1-3, indexed by $k$) alternates between the taxi-bus \eqref{subpiter1} and DV \eqref{subpiter1} subproblems.

\subsection{Discussion and recommendation for practice}

A few extension of the DSC problem \eqref{DSCfeqn1}-\eqref{DSCfeqn6} are discussed below. 

{\bf 1. Sensor collocation.} In applications like air quality sensing, low-cost sensors need to be frequently calibrated through sensor collocation (by staying within the vicinity of a reference monitoring station) while on the move \citep{Saukh2014, Maag2017}. Such a consideration can be reflected in the objective \eqref{DSCfeqn1} by assigning higher weights $w_g$ to grids that contain reference stations. 

\vspace{0.1 in}

{\bf 2. Heterogenous taxi distributions.} In metropolitan areas, taxi services are run by multiple operators, segmenting the service market geographically. In the case of Chengdu city (population: 20 million, area: 3640 km$^2$ with 13 districts), each district such as Longquanyi has their own taxi operators, and the corresponding fleet are only active in local areas. Mathematically, this means the i.i.d. assumption regarding the ball-in-bin model of taxi coverage only applies to a given local operator $l$. To address this issue, we consider the following convex objective:
$$
\max_{n^T_1,\ldots,n^T_L}\sum_{t\in\mathcal{T}} \mu_t \sum_{g\in \mathcal{G}} w_g \left( \sum_{l=1}^Ln^T_l p^{(l)}_{g,t}\right)^{\beta}
$$
\noindent where $n^T_l$ is the size of the instrumented fleet of operator $l$, and $p_{g,t}^{(l)}$ is the binomial parameter associated with operator $l$. In this case, the taxi-related decision variable becomes a vector $(n_1^T,\ldots, n_L^T)$, and the taxi-bus subproblem remains convex. The rest of the DSC formulation incorporating buses and DVs is straightforward.

\section{Computational experiments}\label{secCS}

The drive-by sensing coverage problem with mixed vehicle types is first demonstrated in a real traffic network (Longquanyi District, Chengdu, China) in Sections \ref{subsecCEsetup}-\ref{subsecbeta}. Then, a transferability study is conducted in Section \ref{subsecCEtransf} based on three different urban networks and hundreds of network variations to extrapolate qualitative insights.

\begin{figure}[H]
\centering
\includegraphics[width=.8\textwidth]{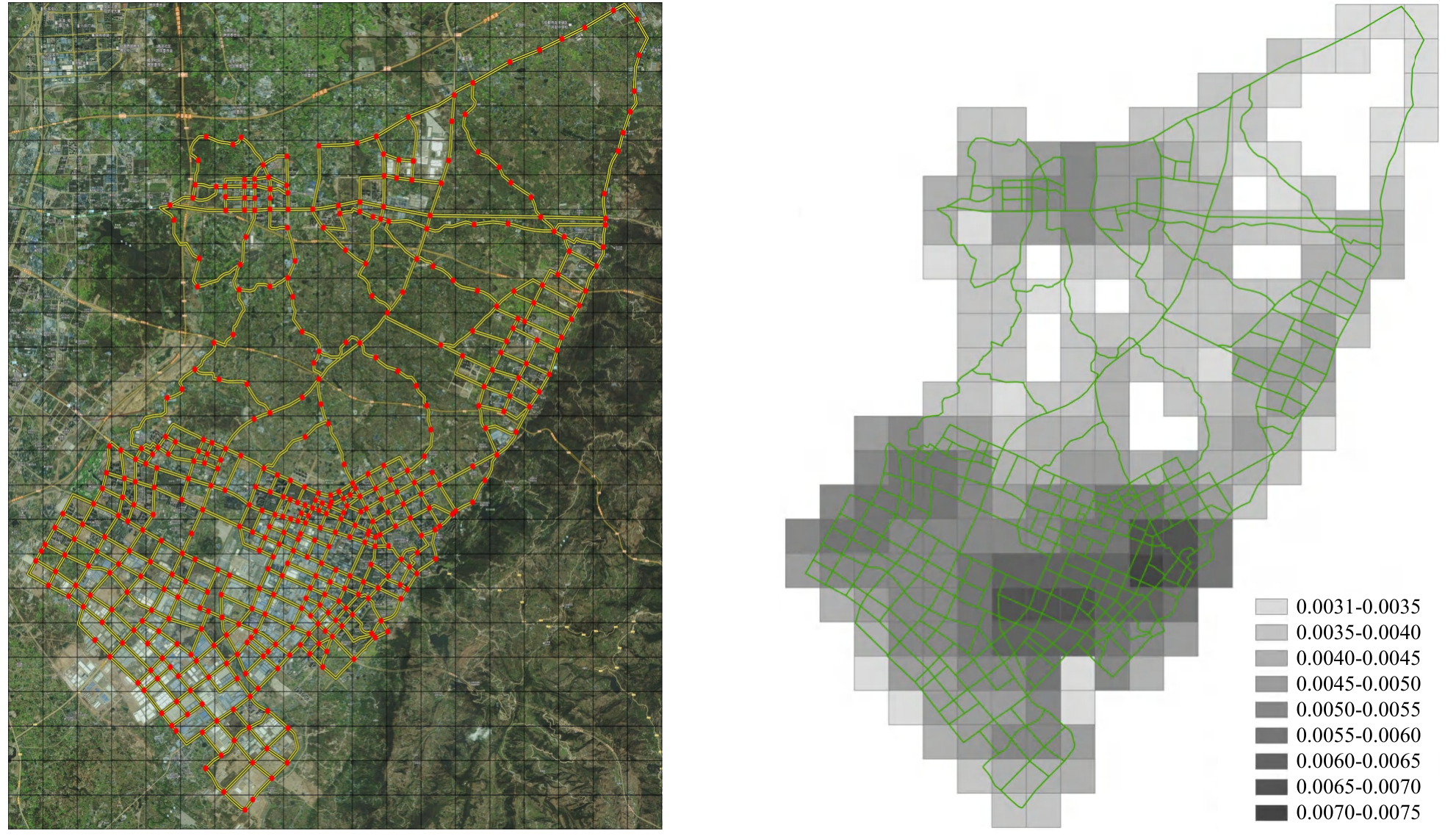}
\caption{Left: Road network and spatial mesh (1km-by-1km) of Longquanyi District, Chengdu. Right: Spatial sensing weights $w_g$ of all the grids covered by the road network.}
\label{fignetwork}
\end{figure}

The Longquanyi area is subject to dense air quality monitoring under a joint industry-government initiative on sustainable urban management. The 252 km$^2$ area is partitioned into grids of size 1km$\times$1km, as shown in Figure \ref{fignetwork} (left). The land use type, population distribution and transport activities vary across the region, constituting a sufficiently sophisticated test case to evaluate the power of mixed-vehicle sensing.

\subsection{Model setup and parameters}\label{subsecCEsetup}

\subsubsection{Spatial and temporal sensing weights}\label{subsubsecstsw}
\begin{figure}[H]
\centering
\includegraphics[width=.6\textwidth]{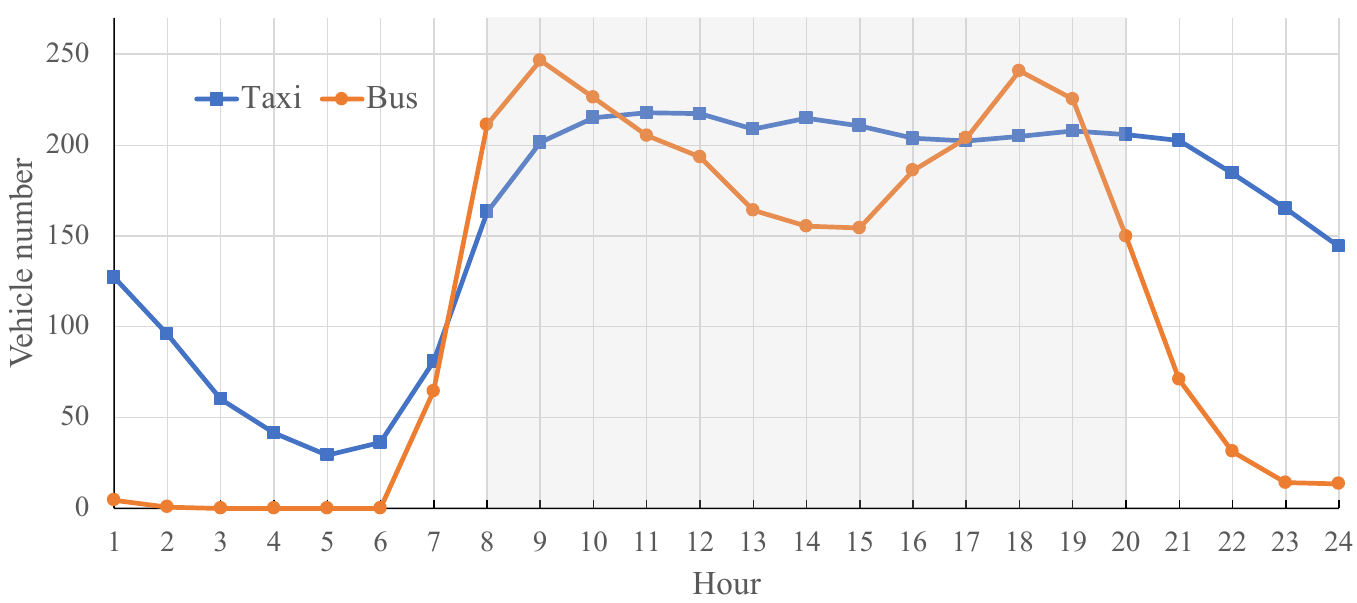}
\caption{Volume of operating taxis and buses in a 24-hour horizon.}
\label{figflow}
\end{figure}

The spatial sensing weights $\{w_g,\,g\in\mathcal{G}\}$ are calculated based on the distribution of five major emission sources (earthwork, factories, traffic, auto services, restaurants), as well as population density. The resulting  weights are shown in Figure \ref{fignetwork} (right). The temporal weights $\mu_t$'s are specified based on two scenarios:
\begin{itemize}
\item {\bf Scenario 1}: A 12-hour horizon (8:00-20:00) with equal weight for each hour ($\mu_t={1\over 12}$, $t=1,\ldots, 12$);
\item {\bf Scenario 2}: A 24-hour horizon with day-time (8:00-20:00) weights quadrupling night-time (21:00-7:00) weights ($\mu_t={1\over 15}$, $t=8,\ldots,20$; $\mu_t={1\over 60}$, otherwise). 
\end{itemize}
\noindent Both scenarios emphasize 8:00-20:00 for air quality monitoring because most emission sources are active during this period. Figure \ref{figflow} shows the time-varying taxi and bus volume in this area. While the bus service is mostly concentrated within 8:00-20:00, the taxi activities remain significant even past midnight. It is expected that taxis will gain more advantage in Scenario 2.

In Sections \ref{subsecORA} and \ref{subsececonms}, the parameter $\beta$ in the utility function $\xi(N)=N^{\beta}$ is set to be $0.2$  following Remark \ref{rmkbeta}. The impact of different values of $\beta$ will be analyzed in Section \ref{subsecbeta}.

\subsubsection{Cost structure}
The mobile air quality sensor (Mode: GLZT\_Y03S) deployed in Longquanyi provides measurements of PM$_{2.5}$, PM$_{10}$, CO, SO$_2$, NO$_2$, and O$_3$ concentrations, at a temporal resolution of 3-5 seconds. These sensors are mounted to taxis, buses and dedicated vehicles. The cost structures are shown in Table \ref{tabcs}.  

\begin{table}[h!]
\centering
  \caption{Unit cost structure for a three-year project period (in CNY). $^{\bf a}$Sensor procurement and maintenance. $^{\bf b}$Compensation for vehicle modification and supporting activities. $^{\bf c}$Vehicle maintenance and energy charges. $^{\bf d}$Vehicle procurement (40k) and insurance (10k). $^{\bf e}$5k/month for 3 years.}
\begin{tabular}{cccccc}
\hline
  & Sensor  & Maintenance & Vehicle cost & Personnel cost & Total (3 years)
  \\\hline
  Taxi & 50k/3 years$^{\bf a}$ & 2k/3 years$^{\bf b}$ & 0 & 0 & 52k
  \\
  Bus & 50k/3 years$^{\bf a}$ & 0 & 0 & 0 & 50k
  \\
  DV & 50k/3 years$^{\bf a}$ & 20k/3 years$^{\bf c}$ & 50k$^{\bf d}$ & 180k/3 years$^{\bf e}$ & 300k
  \\\hline
  \end{tabular}
  \label{tabcs}
  \end{table}

\subsubsection{Operational parameters of sensing vehicles}\label{subsubsecOP}

{\bf Taxis:} A total of 264 taxis are registered with the local operator. Their GPS trajectories in March 2021 are used to fit and validate the linear approximation \eqref{MT}, as detailed in Section \ref{subsecAppT}. 

\vspace{0.1in}

 \noindent {\bf Buses:} A total of 56 bus lines are considered for the case study, whose routes are shown in Figure \ref{figTaxiBus}(b). The GPS data of 253 buses are used to estimate the average en route service times $T_j^s(t)$, turn-around times $T_j^a(t)$, and the number of operating buses $\lambda_j(t)$ per line. 
 
\vspace{0.1in}

\noindent {\bf Dedicated vehicles:} The operating hours of DVs are set to be 8:00-20:00. The speed of the DVs is set to be 30 km/h, but it is straightforward to incorporate spatially and temporally heterogenous speed distributions.

\subsection{Evaluation of the optimization results}\label{subsecORA}

Following Section \ref{subsecQSTC}, the goodness-of-coverage is quantified in two ways: (1) the {\it space-time weighted sensing utility} (STWSU) $\Phi$, which is the optimization objective of the DSC problems; and (2) the {\it KL-divergence} of the actual distribution from the target distribution of the sensing utilities $\xi(N_{g,t})$. 

\begin{figure}[H]
\centering
\includegraphics[width=.85\textwidth]{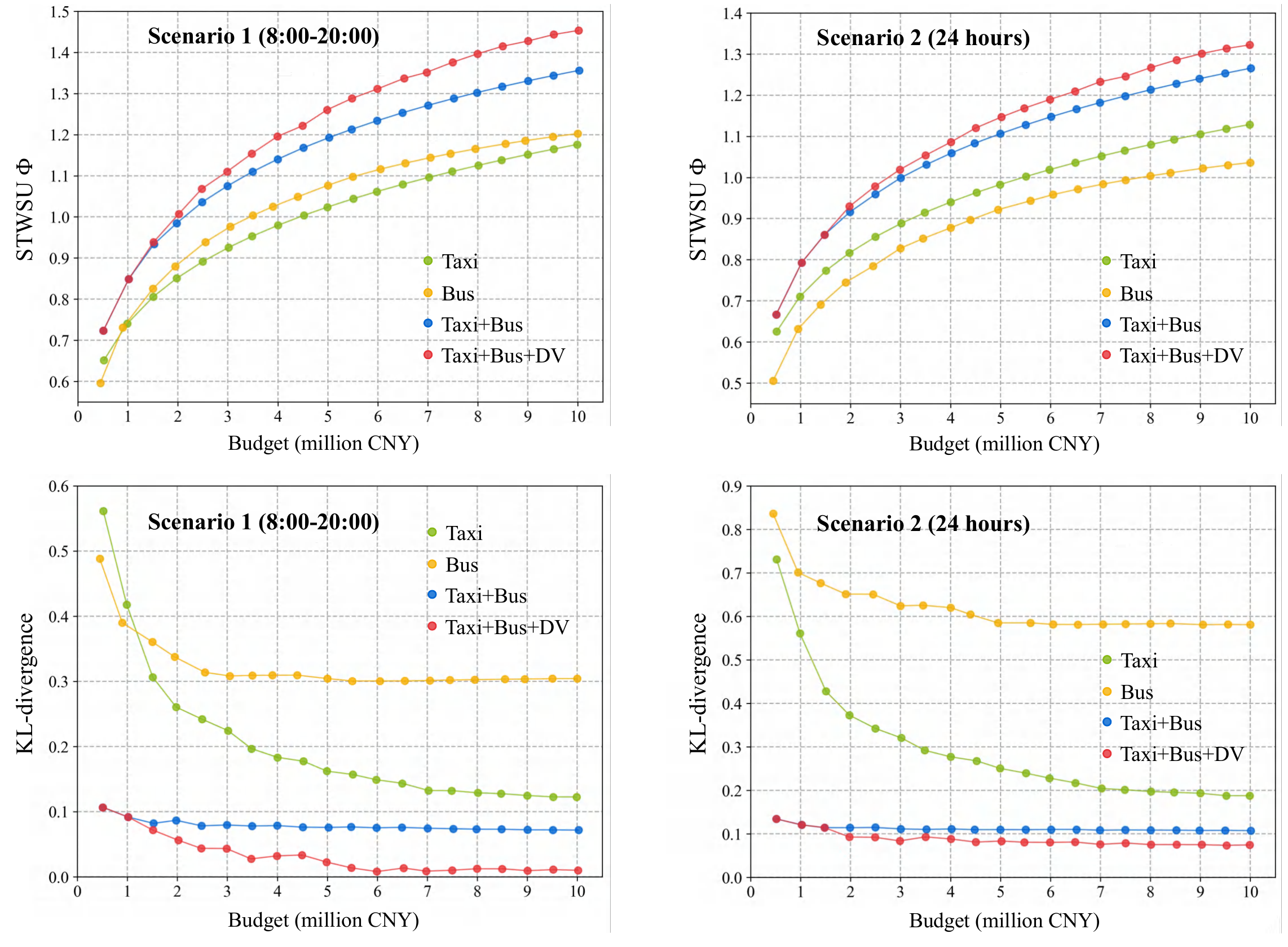}
\caption{Goodness-of-coverage of four fleet combinations under various budgets. First row: STWSU $\Phi$; second row: KL-divergence of the actual distribution from the target distribution of utilities.}
\label{figcurves}
\end{figure}

Four vehicle fleet combinations are considered for the DSC problem: (1) taxi; (2) bus; (3) taxi+bus; and (4) taxi+bus+DV. 
As Figure \ref{figcurves} shows, under both scenarios 1 and 2, mixed fleets significantly outperform single fleets in terms of both $\Phi$ and KL-div. With the same budget, the STWSU of the taxi+bus+DV mode exceeds that of single modes, by 11.1\%-20.9\% in Scenario 1 and 6.7\%-17.7\% in Scenario 2. Moreover, the sensing powers of buses and DVs are lower in Scenario 2, as it includes night time (21:00-7:00) during which they are inactive compared to taxis. 

 It is apparent that single-mode fleets are inadequate to approximate target distributions for the lack of flexibility in spatial coverage: Taxis follow a priori distributions and buses collect utilities only along their routes. It is remarkable, however, that the two combined drastically reduces the KL-div., even with very low budgets. This suggests a high level of synergy between taxis and buses achieved through the optimization procedure. Furthermore, the addition of DVs can bring the KL-div. to a very low value in Scenario 1, while they are inactive during night times in Scenario 2. 


To achieve intuitive visualization of the spatial coverage, we define the {\it time-weighted sensing utility} (TWSU) and {\it time-averaged gap} (TAG):
\begin{equation}\label{TWSUTAG}
\hbox{TWSU:}~\sum_{t\in\mathcal{T}}\mu_t \xi(N_{g,t}),\quad \hbox{TAG:}~{1\over  |\mathcal{T}|}\sum_{t\in\mathcal{T}} {\hbox{AD}_{g,t}-\hbox{TD}_{g,t}\over \hbox{TD}_{g,t}}\times 100\%,\qquad \forall g\in\mathcal{G}
\end{equation}
\noindent where $\hbox{AD}_{g,t}$ and $\hbox{TD}_{g,t}$ are actual and target distributions given by \eqref{ADTD}. The former visualizes the spatial distribution of sensing utilities and the latter shows the relative discrepancy between the two distributions.

\begin{figure}[H]
\centering
\includegraphics[width=\textwidth]{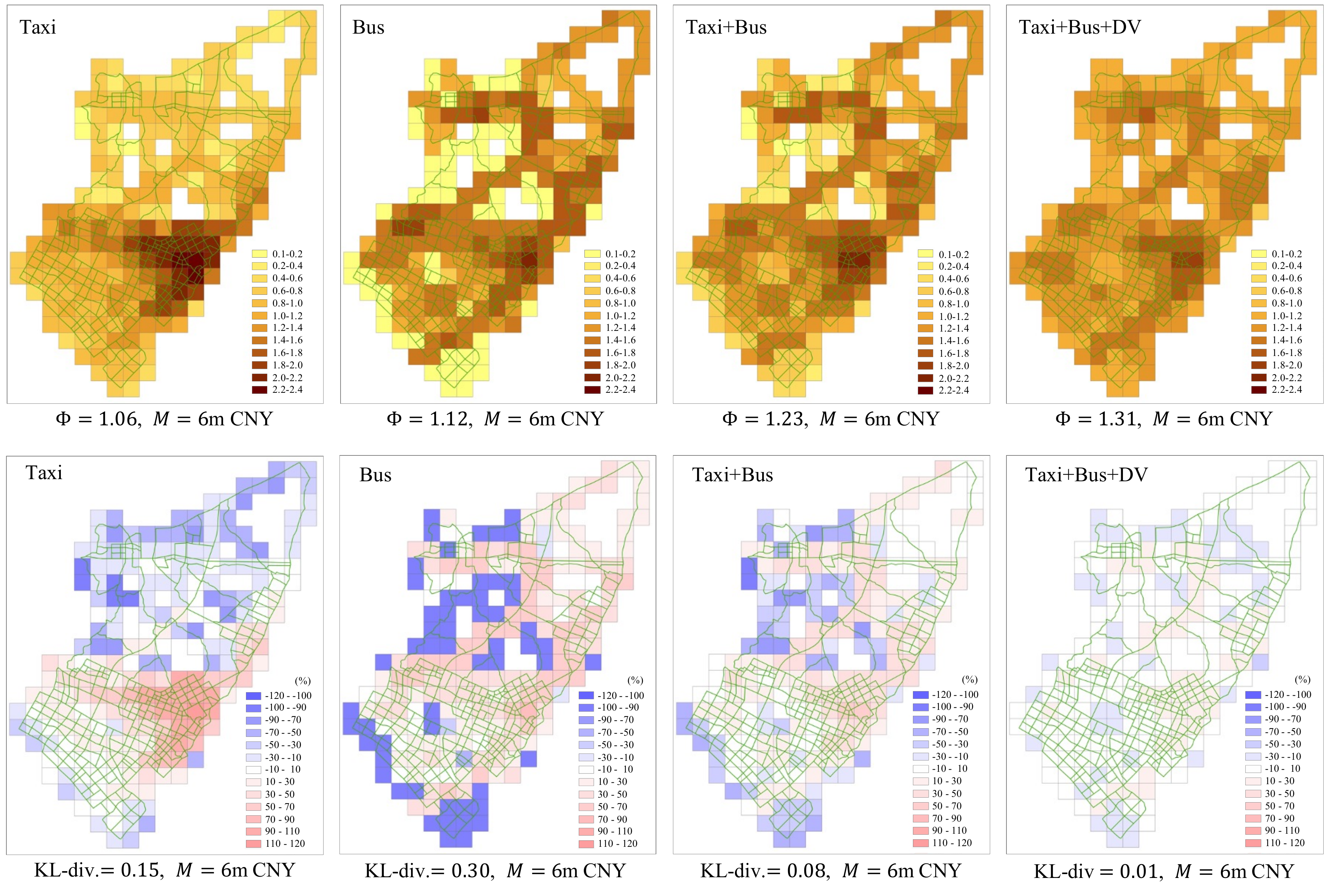}
\caption{Visualization of TWSU (first row) and TAG (second row) from \eqref{TWSUTAG} in Scenario 1 with the same budget of 6m CNY.}
\label{figgc}
\end{figure}

Figure \ref{figgc} shows the TWSU and TAG for four vehicle combinations in Scenario 1 under the same budget of 6m CNY. Taxis provide an unbalanced spatial coverage, concentrating primarily on densely populated areas. Bus coverage has a large spatial extent, but leaves a number of grids unvisited, leading to locally unbalanced coverage. The taxi and bus fleets complement each other with more balanced distribution and only a few spots less covered. These spots are further addressed by DVs following the designed routes, rendering a high level of approximation to the target distribution.

\subsection{The economies of mixed sensing}\label{subsececonms}

\begin{figure}[H]
\centering
\includegraphics[width=\textwidth]{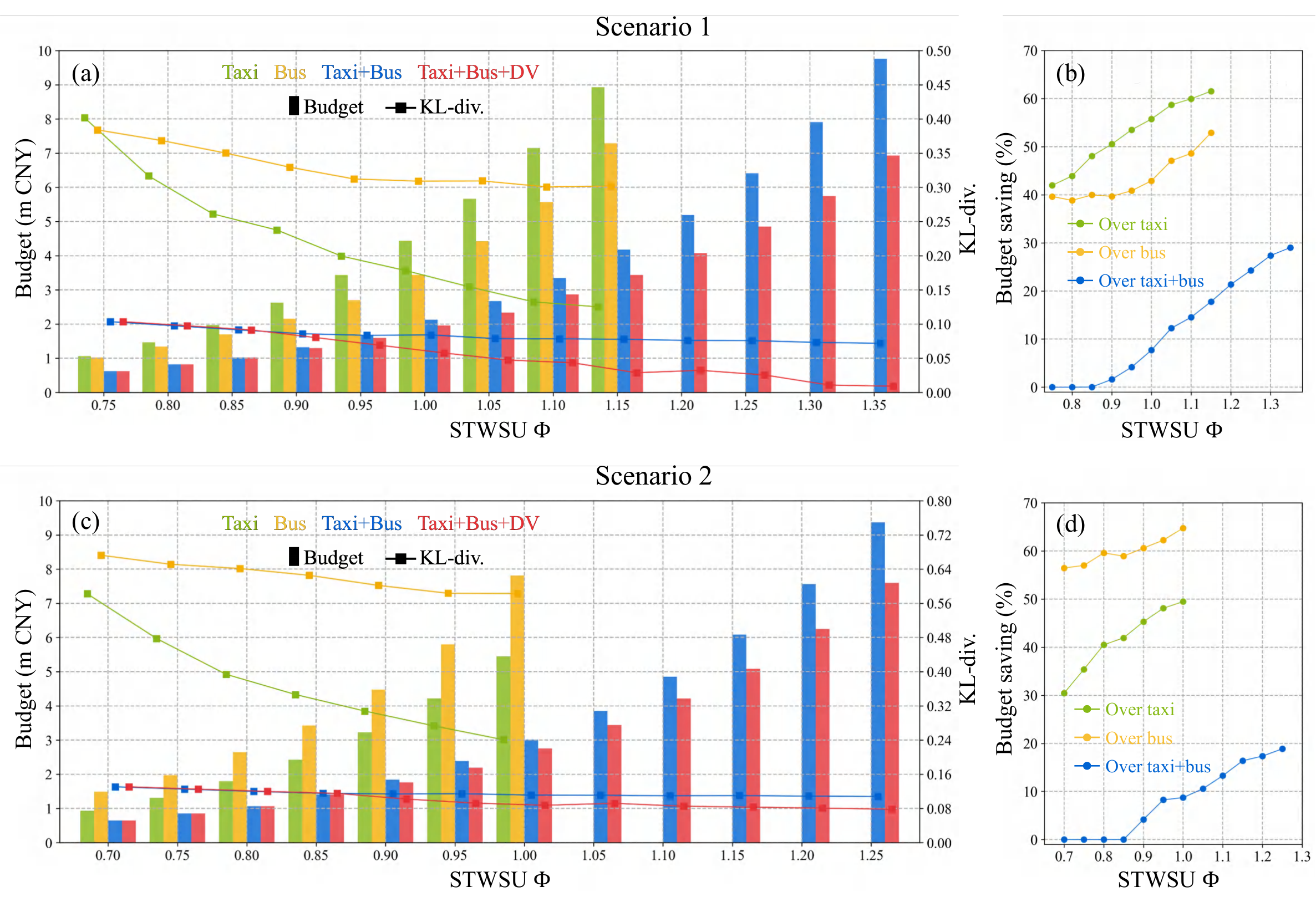}
\caption{(a) \& (c): Budget and KL-div. of four vehicle combinations under the same STWSU $\Phi$. (b) \& (d): Budget saving (\%) of the Taxi+bus+DV mode over Taxi, Bus, and Taxi+bus modes.}
\label{figsavings}
\end{figure}

The economies of mixed sensing refers to the budget savings achievable while maintaining the same or even higher sensing performance. Figure \ref{figsavings}(a) \& (c) show the budgets needed by the four fleet combinations to achieve certain values of $\Phi$ (bar chart), as well as the corresponding KL-div. (lines). The budgets of mixed sensing (taxi+bus or taxi+bus+DV) are considerably lower than single modes, with savings ranging from 39\%-61\% in Scenario 1 (subfigure b) and 30\%-65\% in Scenario 2 (subfigure d). Such savings grow with the target STWSU $\Phi$. It is also interesting that considering DVs can achieve considerable sensing economies (compared to taxi+bus), despite the fact that its unit cost is much higher.

\begin{figure}[H]
\centering
\includegraphics[width=\textwidth]{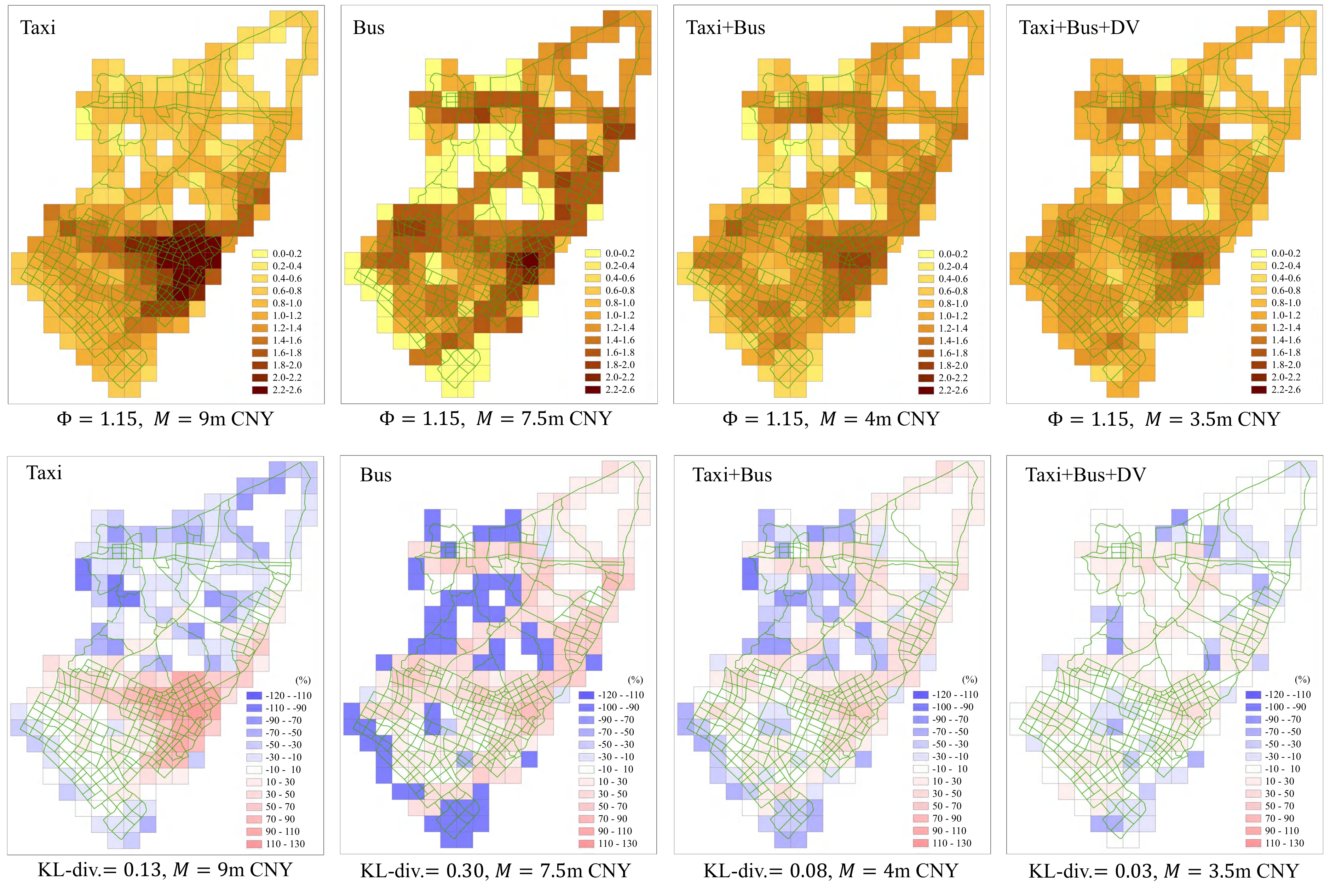}
\caption{Visualization of TWSU (first row) and TAG (second row) from \eqref{TWSUTAG} under the same STWSU $\Phi=1.15$, in Scenario 1.}
\label{figeconmg}
\end{figure}

Another important observation from Figure \ref{figsavings} is that, for the same STWSU $\Phi$, the KL-div. is always smaller when there are more fleet types in the mix. This suggests that the diversified mobility patterns of the fleet can be explored to achieve far greater sensing power than single modes. For intuitive visualization, Figure \ref{figeconmg} shows the spatial distribution of utilities (TWSU) and distribution discrepancies (TAG) when the four fleet combinations achieve the same $\Phi=1.15$ under Scenario 1. The taxi-bus-DV mode achieves the best sensing performance in terms of sensing distribution, with less than half the budget of taxi or bus.

\subsection{Sensitivity analysis of $\beta$} \label{subsecbeta}

\begin{figure}[H]
\centering
\includegraphics[width=\textwidth]{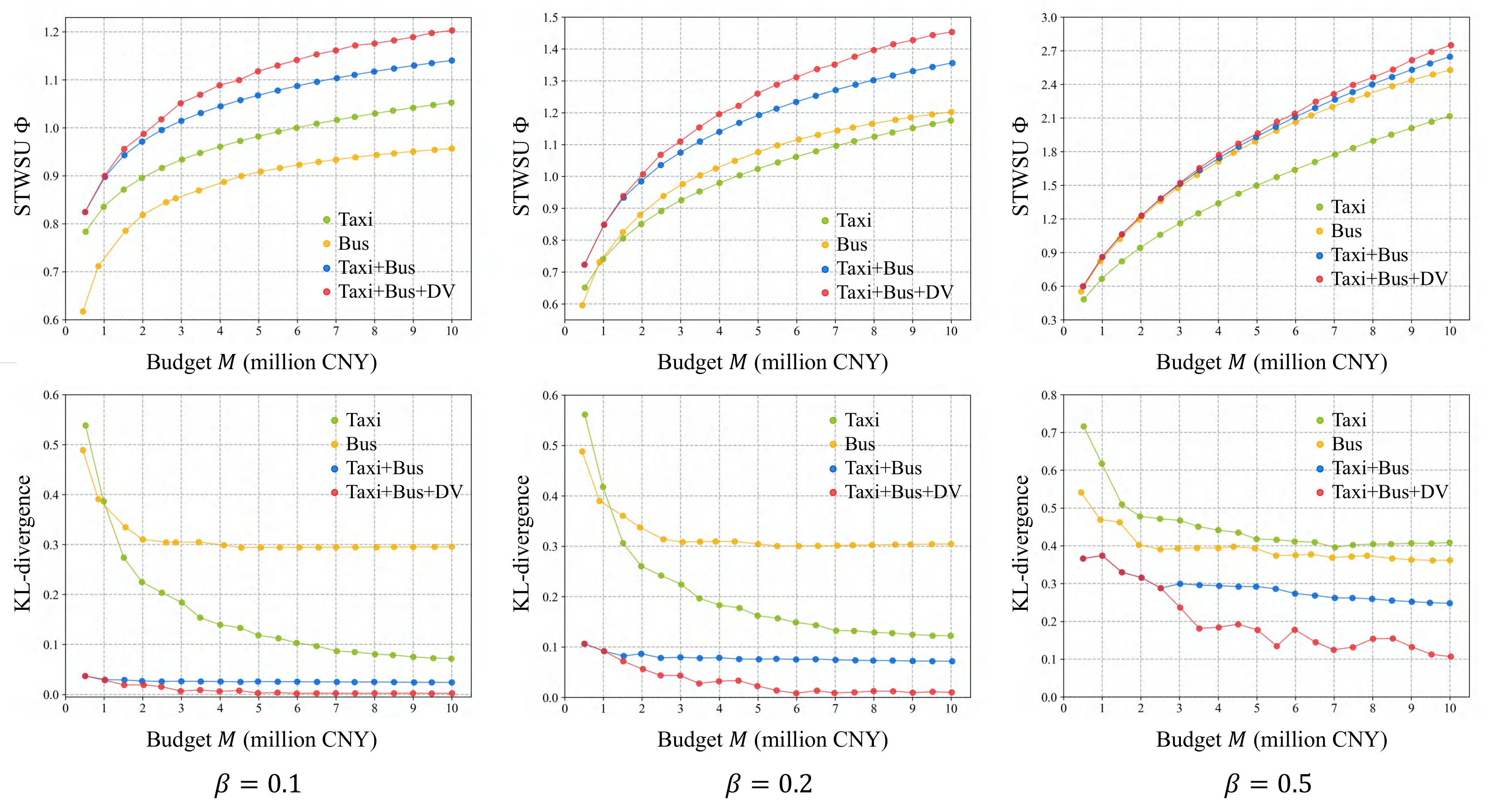}
\caption{Performance measures (STWSU and KL-div.) with different $\beta$ values.}
\label{figbeta}
\end{figure}

The results presented so far are all based on the sensing utility function $\xi(N)=N^{\beta}$ where $\beta=0.2$. To understand the impact of $\beta$ on the performance of different vehicle fleet, in Figure \ref{figbeta} we plot STWSU and KL-divergence with $\beta=0.1, 0.2, 0.5$, under Scenario 1. Recall from Eqn \ref{proptoeqn} that a small $\beta$ promotes relatively uniform distribution of utilities $\xi(N_{g,t})$, while for a large $\beta$ the utilities are relatively concentrated in high-weighting grids. From the standpoint of diminishing marginal utility in $\xi(N)=N^\beta$, a small $\beta$ discourages repeated sensing coverage or requires low sensing frequency.

Firstly, all the $M$-$\Phi$ curves are concave and, in fact, follow the relationship $\Phi=aM^{\beta}+b$ for some $a, b>0$. Secondly, in terms of both STWSU and KL-div., the advantage of mixed sensing is more pronounced for smaller $\beta$, which requires relatively uniform coverage. Thirdly, recall from Section \ref{subsubsecTD} that the KL-div. indicates the loss of sensing efficacy caused by the fleet's mobility constraints. For all values of $\beta$, the mixed fleets (Taxi+Bus, Taxi+Bus+DV) yield satisfactory approximation to the target distribution, and their KL-div. stays relatively low even under small budgets. Finally, comparing Taxi+Bus with Taxi+Bus+DV reveals that the DVs are used in the optimal solution for $M\geq 1.0$m CNY ($\beta=0.1$), $M\geq 1.5$m CNY ($\beta=0.2$) and $M\geq 3.0$m CNY ($\beta=0.5$). This suggests that for small $\beta$, the DVs are essential to ensure balanced coverage, and are required even with low budgets (e.g. 1.5 m).



\subsection{Transferability study}\label{subsecCEtransf}

\begin{figure}[H]
\centering
\includegraphics[width=\textwidth]{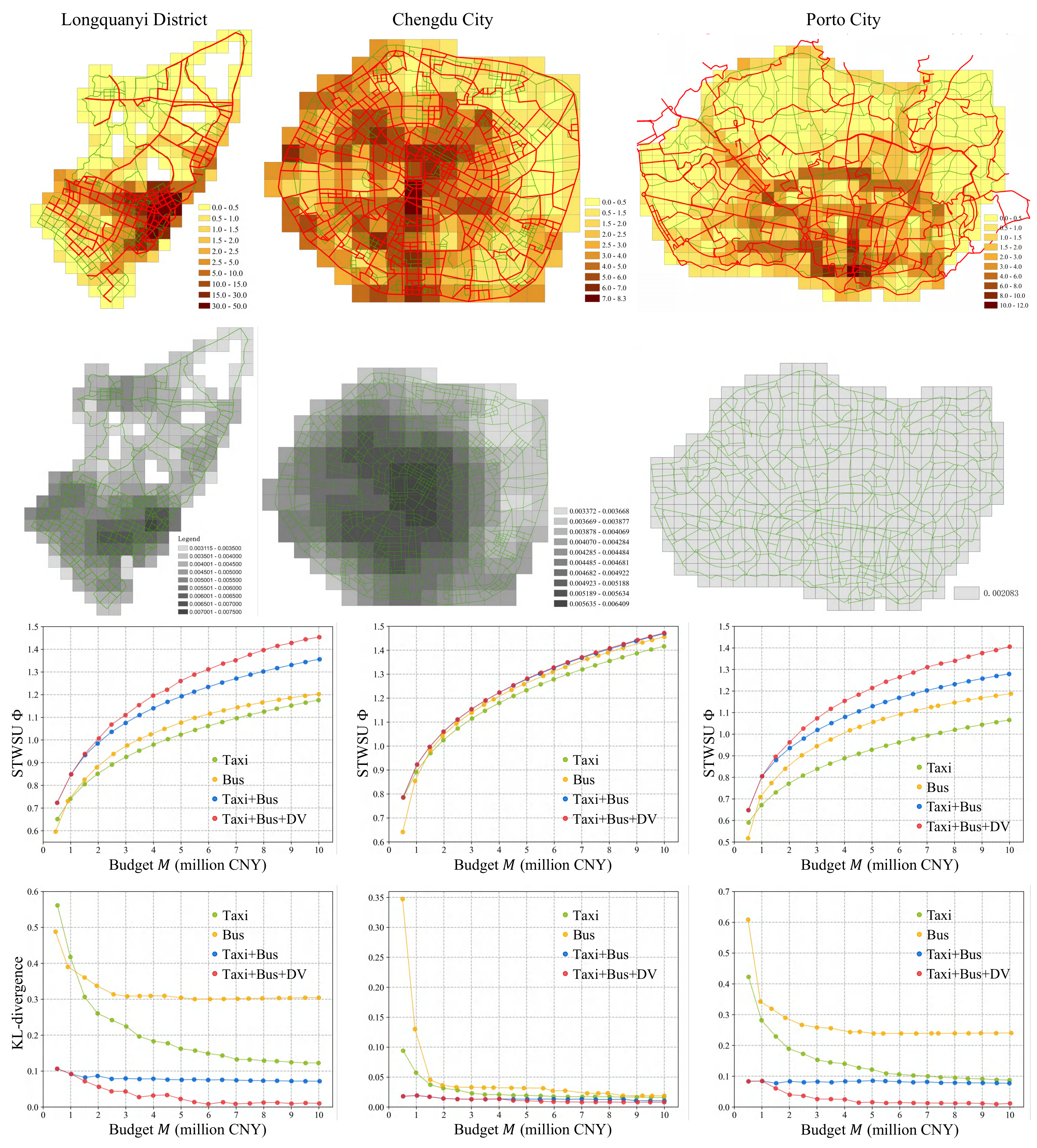}
\caption{First row: The three study areas where green lines are road networks, red lines are bus routes, and colored grids represent mean taxi counts (with fleet size $n^T=80$) during 8-9am. Second row: Spatial sensing weights $w_g$'s. Third \& forth rows: The STWSU $\Phi$ and KL-div. by solving the DSC problems with different budgets.}
\label{fig3cities}
\end{figure}

To further understand the sensing power of mixed vehicle fleets beyond the context of a fixed network configuration, we consider two more real-world networks (Figure \ref{fig3cities}): The cities of Chengdu (China) and Porto (Portugal). We use taxi GPS data to fit the expectation $N_{g,t}^T$ of the binomial distribution \eqref{MT}, as well as bus operational data to estimate the expected coverage $N_{g,t}^B$ \eqref{MB}. The spatial sensing weights for Chengdu (based on a similar approach as Longquanyi) and Porto (uniform weights for the lack of relevant data) are also shown in Figure \ref{fig3cities}. 

The DSC problems with mixed fleets are solved for the three networks, and the STWSU $\Phi$ and KL-div. are shown in Figure \ref{fig3cities}. The results for Longquanyi and Porto are quite similar, as both networks are characterized by highly unbalanced taxi distributions and numerous spatial grids not covered by bus routes. In contrast, both taxi and bus route distributions in Chengdu offer adequate coverage of the target area, resulting in marginal improvement of mixed fleets.

To gain further analytical insights into mixed-fleet sensing, we generate a large number of network configurations (samples) by randomly removing bus routes, thereby using the bus network coverage as the control variable to enrich our datasets\footnote{Ideally we would use both bus networks and taxi distributions, but the latter are more difficult to populate.}. For a given network configuration, we define:
\begin{itemize}
\item $W(T)$: Sum of sensing weights ($w_g$'s) of grids that have a taxi coverage within the upper $60\%$ among all; 
\item $W(B)$: Sum of sensing weights ($w_g$'s) of grids that intersect any bus route;
\item $W(T+B)$:  Sum of sensing weights ($w_g$'s) of grids that either intersect bus routes, or have a taxi coverage within the upper $60\%$ among all;
\item $W(T+B+D)$: Sum of sensing weights $w_g$'s of grids that can be covered by either taxis, buses or DVs. By default, $W(T+B+D)$ is set to be 1. 
\end{itemize}

\noindent Here, $W(T)$, $W(B)$ and $W(T+B)\in (0,\,1)$ are simple and easy-to-calculate \underline{s}patial \underline{c}overage \underline{i}ndicators (SCIs), prior to any sensor allocation procedures. The following relative SCI differences are calculated.
\begin{equation}\label{rSCId}
{W(T+B)-W(T)\over W(T)},~~{W(T+B)-W(B)\over W(B)}, ~~ {W(T+B+D)-W(T+B)\over W(T+B)}
\end{equation}

Next, for each network configuration, the DSC problem is solved with various fleet combinations to obtain 
\begin{equation}\label{rSG}
{\Phi(T+B)-\Phi(T)\over \Phi(T)},~~{\Phi(T+B)-\Phi(B)\over \Phi(B)},~~ {\Phi(T+B+D)-\Phi(T+B)\over \Phi(T+B)},
\end{equation}
\noindent which represent the relative sensing gains (STWSU), with obvious meaning of notations. 

\begin{figure}[H]
\centering
\includegraphics[width=\textwidth]{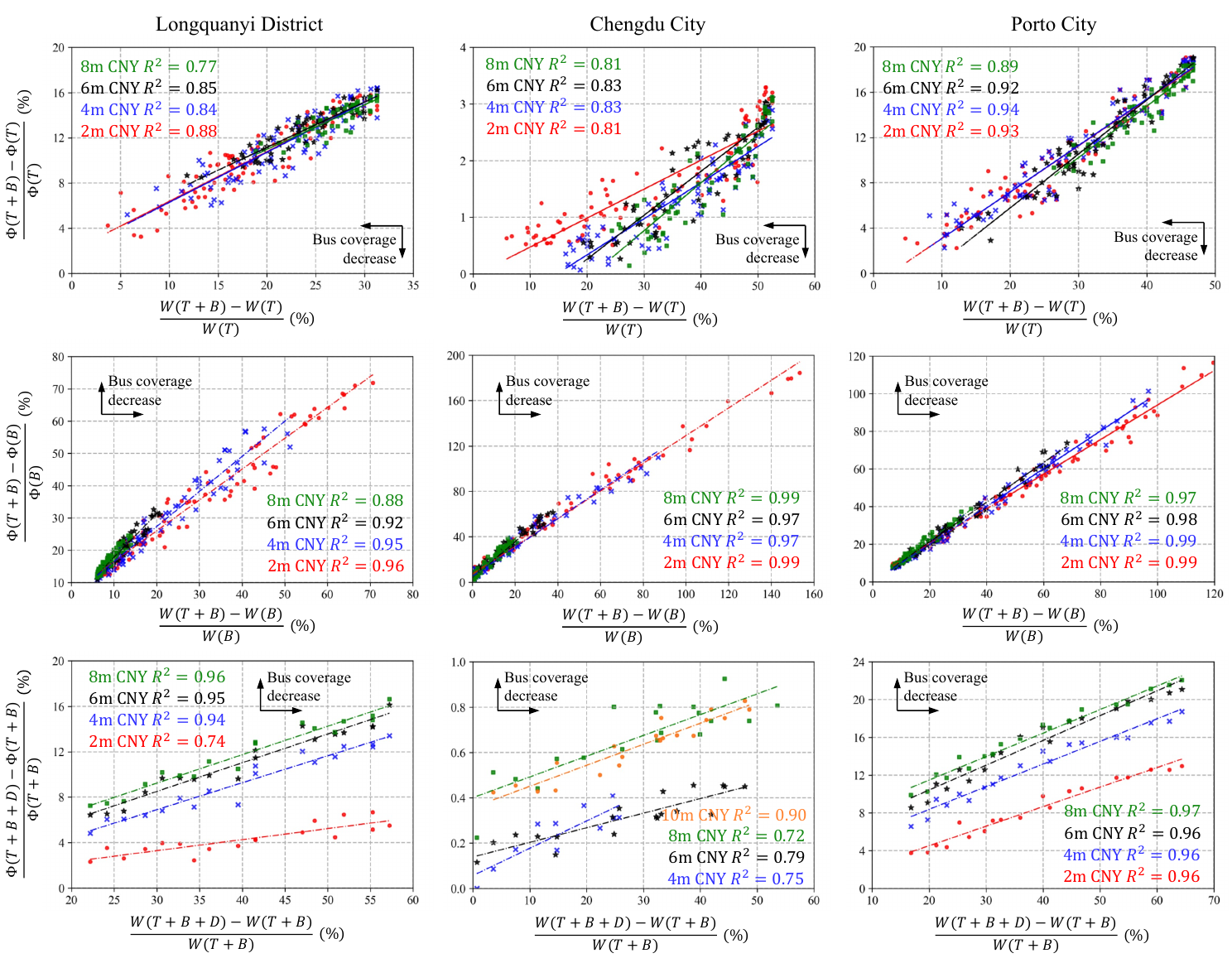}
\caption{Regression analyses between relative SCI differences \eqref{rSCId} and relative sensing gains \eqref{rSG}. The scatter plots are based on samples populated by randomly removing bus routes from the original networks. }
\label{figregression}
\end{figure}

The relative SCI differences \eqref{rSCId} and relative sensing gains \eqref{rSG} are plotted in Figure \ref{figregression}. We observe that they display linear relationships (with high R-squares) at various budget levels for all three networks. 

\begin{itemize}
\item First row: As bus network coverage $W(B)$ gradually decreases, both $\big(W(T+B)-W(T)\big)/W(T)$ and $\big(\Phi(T+B)-\Phi(T)\big)/ \Phi(T)$ decrease, but the latter is within a relatively narrow range $(0\%,\,4\%)$ in the case of Chengdu. This means that the spatial extent of bus networks play an important role in complementing the sensing power of taxi fleets, unless the taxi fleets already provide good coverage (e.g. Chengdu). 

\item Second row: As the level of bus network coverage $W(B)$ gradually decreases, both $\big(W(T+B)-W(B)\big)/ W(B)$ and $\big(\Phi(T+B)-\Phi(B)\big)/ \Phi(B)$ increase. The latter is mostly within $(0\%,\,100\%)$. This means taxi fleets could substantially boost the sensing power of bus fleets, especially when the bus networks are sparse. 

\item Third row: As the level of bus network coverage $W(B)$ gradually decreases, both $\big(W(T+B+D)-W(T+B)\big)/ W(T+B)$ and $\big(\Phi(T+B+D)-\Phi(T+B)\big)/\Phi(T+B)$ increase, but the latter is within a relatively narrow range $(0\%,\,1\%)$ in the case of Chengdu. This means that DVs can significantly improve the sensing utility (by up to 16\% in Longquanyi and 22\% in Porto), unless the taxi distributions are already adequate (Chengdu); in addition, their contributions are more pronounced under more sparse bus networks and higher budgets. 
\end{itemize}

\section{Discussion and conclusion}\label{secconclude}

The sensing power of mixed fleet (taxis, buses, DVs) is explored via an optimization approach termed the drive-by sensing coverage (DSC) problem. The sensing efficacy is assessed using space-time weighted sensing utility $\Phi$, as well as gap (KL-divergence) to a target distribution of sensing coverage. The distinct mobility characteristics of taxis, buses and DVs poses significant challenge to the DSC problem, while at the same time showing considerable synergistic potential in offering balanced and targeted sensing capabilities. 

\subsection{Qualitative insights}
The following are drawn from our various experiments and presented in a manner that emphasizes insights beyond numerical results. 

\begin{itemize}
\item In most cases, buses collect more information (higher sensing utilities $\Phi$) than taxis for the greater spatial extent of bus routes, but bus-based sensing could be quite biased (large KL-div.) against areas with no bus routes. Combining both offers satisfactory sensing performance in terms of both quantity and balance (fairness). 

\item Using a mix of vehicle types can easily achieve good approximation to a target distribution, even with low budgets and low sensing utilities $\Phi$. Thus, for large areas with diversified sensing needs, the simultaneous use of taxis and buses (if not DVs) is highly recommended. 

\item In most cases, the spatial extent and density of bus networks are strongly and positively correlated to the sensing power of mixed fleets, unless the taxi distributions already offer sufficient coverage. For areas with highly developed bus networks, the consideration of buses in the mix is highly beneficial. 

\item Although DV is a much more expensive alternative, it is proven to bring considerable economies to drive-by sensing. It is particularly important for locations insufficiently covered by taxis or bus routes, and for sensing applications with low temporal frequency requirements (small $\beta$).

\item The time window and temporal sensing weights also play a role in the choice of mixed fleets, which were neglected in existing studies. In this case, the operational hours of various fleets should be a key factor in the problem. 
\end{itemize}

\subsection{Other variants of the problem and future research}\label{secmore}

\begin{itemize}

\item {\bf Boosting the sensing power of taxis.} One cost-effective way to increase the spatial coverage of taxis is through mechanism design on ridesourcing platforms, including taxi-rider matching and pricing. The former is essential to address the spatial coverage of instrumented fleet, and the latter addresses player incentives and budget balance.

\item {\bf Flexible bus-line binding.} For cases with low sensing frequency requirements (low $\beta$), an instrumented bus can be relocated to a different line at the terminal, such that a single sensor can be used to serve multiple lines, covering a much larger area over a longer period of time \citep{DH2023}.

\item {\bf Area assignment for DVs.} In practice, dedicated vehicles need to perform a spectrum of tasks beyond drive-by sensing, including unscheduled on-site investigations. The fixed routes and itineraries derived in this work should be replaced with a more flexible area assignment module, with a probabilistic account on the space-time coverage within each area.

\end{itemize}

\section*{Acknowledgement}
This work is supported by the National Natural Science Foundation of China through grant 72071163, and the Natural Science Foundation of Sichuan Province through grant 2022NSFSC0474.


\begin{thebibliography}{99}

\bibitem[Ali and Dyo et al., 2017]{AD2017} Ali, J., Dyo, V., 2017. Coverage and Mobile Sensor Placement for Vehicles on Predetermined Routes: A Greedy Heuristic Approach:, in: Proceedings of the 14th International Joint Conference on E-Business and Telecommunications. Presented at the 14th International Conference on Wireless Networks and Mobile Systems, SCITEPRESS - Science and Technology Publications, Madrid, Spain, pp. 83–88. 

\bibitem[Alsina-Pag\'es et al., 2017]{Alsina-Pages2017} Alsina-Pag\'es, R., Hernandez-Jayo, U., Alías, F., Angulo, I., 2017. Design of a Mobile Low-Cost Sensor Network Using Urban Buses for Real-Time Ubiquitous Noise Monitoring. Sensors 17(1),57.

\bibitem[Anjomshoaa et al., 2018]{ADRMdR2018} Anjomshoaa, A., Duarte, F. Rennings, D., Matarazzo, T.J., deSouza, P., Ratti, C. 2018. City Scanner: Building and Scheduling a Mobile Sensing Platform for Smart City Services. IEEE Internet of Things Journal, 5(6), 4567-4579.

\bibitem[Asprone et al., 2021]{ADFS2021} Asprone, D., Di Martino, S., Festa, P., Starace, L. L. L.,2021. Vehicular crowd-sensing: a parametric routing algorithm to increase spatio-temporal road network coverage. International Journal of Geographical Information Science, 35(9), 1876-1904.

\bibitem[Chen et al., 2020]{C2020} Chen, X., Xu, S., Han, J., Fu, H., Pi, X., Joe-Wong, C., Li, Y., Zhang, L., Noh, H.Y., Zhang, P., 2020. PAS: Prediction-Based Actuation System for City-Scale Ridesharing Vehicular Mobile Crowdsensing. IEEE Internet of Things Journal 7 (5), 3719-3734.

\bibitem[Chen et al., 2017]{CLGZX2017} Chen, Y., Lv, P., Guo, D., Zhou, T., Xu, M., 2017. Trajectory segment selection with limited budget in mobile crowd sensing. Pervasive and Mobile Computing 40, 123-138. 

\bibitem[Cruz et al., 2020]{Cruz2020} Cruz, P., Couto, R.S., Costa, L.H.M.K., Fladenmuller, A., Dias de Amorim, M., 2020. A delay-aware coverage metric for bus-based sensor networks. Computer Communications 156, 192-200.

\bibitem[Cruz et al., 2020]{CCCFD2020} Cruz, P., Couto, R.S., Costa, L.H.M.K., Fladenmuller, A., de Amorim, M.D., 2020. Per-Vehicle Coverage in a Bus-Based General-Purpose Sensor Network. IEEE Wireless Communications Letters 9 (7), 1019-1022.

\bibitem[Cruz Caminha et al., 2018]{CDMFD2018} Cruz Caminha, P., de Souza Couto, R., Maciel Kosmalski Costa, L., Fladenmuller, A., Dias de Amorim, M., 2018. On the coverage of bus-based mobile sensing. Sensors 18 (6), 1976. 

\bibitem[Dai and Han, 2023]{DH2023} Dai, Z, Han, K, 2023. Exploring the drive-by sensing power of bus fleet through active scheduling. Transportation Research Part E: Logistics and Transport Review, 171, 103029.

\bibitem[Dang et al., 2013]{DGM2013} Dang, D.-C., Guibadj, R.N., Moukrim, A., 2013. An effective PSO-inspired algorithm for the team orienteering problem. European Journal of Operational Research 229, 332–344.



\bibitem[deSouza et al., 2020]{dADKKR2020} deSouza, P., Anjomshoaa, A., Duarte, F., Kahn, R., Kumar, P., Ratti, C., 2020. Air quality monitoring using mobile low-cost sensors mounted on trash-trucks: Methods development and lessons learned. Sustainable Cities and Society 60, 102239.

\bibitem[Du et al., 2015]{DCYLGS2015} Du, R., Chen, C., Yang, B., Lu, N., Guan, X., Shen, X.,2015. Effective Urban Traffic Monitoring by Vehicular Sensor Networks. IEEE Transactions on Vehicular Technology, 64(1), 273-286.

\bibitem[Du et al., 2019]{DSXVF2019} Du, R., Santi, P., Xiao, M., Vasilakos, A.V., Fischione, C., 2019. The Sensable City: A Survey on the Deployment and Management for Smart City Monitoring. IEEE Communications Surveys \& Tutorials 21 (2), 1533-1560.

\bibitem[Eriksson et al., 2008]{Eriksson2008} Eriksson, J., Girod, L., Hull, B., Newton, R., Madden, S., Balakrishnan, H., 2008. The pothole patrol: using a mobile sensor network for road surface monitoring. in: Proceeding of the 6th international conference on Mobile systems, applications, and services - MobiSys ’08, New York, NY, USA, 29-39.

\bibitem[Fan et al., 2021]{FJLQG2021} Fan, G., Jin, H., Liu, Q., Qin, W., Gan, X., Long, H., Fu, L., Wang, X., 2021a. Joint Scheduling and Incentive Mechanism for Spatio-Temporal Vehicular Crowd Sensing. IEEE Transactions on Mobile Computing 20 (4), 1449-1464.

\bibitem[Fan et al., 2021]{FZGJGW2021}  Fan, G., Zhao, Y., Guo, Z., Jin, H., Gan, X., Wang, X., 2021b. Towards Fine-Grained Spatio-Temporal Coverage for Vehicular Urban Sensing Systems. in: IEEE INFOCOM 2021 - IEEE Conference on Computer Communications, IEEE, Vancouver, BC, Canada, 1-10.

\bibitem[Gao et al., 2016]{Gao2016} Gao, Y., Dong, W., Guo, K., Liu, Xue, Chen, Y., Liu, Xiaojin, Bu, J., Chen, C., 2016. Mosaic: A low-cost mobile sensing system for urban air quality monitoring. in: IEEE INFOCOM 2016 - The 35th Annual IEEE International Conference on Computer Communications, San Francisco, CA, USA, 1-9. 


\bibitem[Glock and Meyer, 2020]{GM2020} Glock, K., Meyer, A., 2020. Mission Planning for Emergency Rapid Mapping with Drones. Transportation Science, 54, 534-560. 



\bibitem[He et al., 2015]{HCL2015} He, Z., Cao, J., Liu, X., 2015. High quality participant recruitment in vehicle-based crowdsourcing using predictable mobility. in: 2015 IEEE Conference on Computer Communications (INFOCOM), Kowloon, Hong Kong, 2542-2550.

\bibitem[Ji et al., 2023a]{JHG2023} Ji, W., Han, K., Ge, Q., 2023. Extended team orienteering problem: Algorithms and applications, arXiv: 2307.02397.

\bibitem[Ji et al., 2023b]{JHL2023} Ji, W., Han, K., Liu, T., 2023. A survey of urban drive-by sensing: An optimization perspective. Sustainable Cities and Society, 99, 104874.

\bibitem[Ji et al., 2023c]{JHL2023partC} Ji, W., Han, K., Liu, T., 2023. Trip-based mobile sensor deployment for drive-by sensing with bus fleets. Transportation Research Part C: Emerging Technologies, 157, 104404.

\bibitem[Kaivonen and Ngai, 2020]{Kaivonen2020} Kaivonen, S., Ngai, E.C.-H., 2020. Real-time air pollution monitoring with sensors on city bus. Digital Communications and Networks 6 (1), 23-30.

\bibitem[Ke et al., 2015]{KZLC2015} Ke, L., Zhai, L., Li, J., Chan, F.T.S., 2015. Pareto mimic algorithm: An approach to the team orienteering problem. Omega 61, 155-166.

\bibitem[Khan et al., 2016]{KGB2016} Khan, J.A., Ghamri-Doudane, Y., Botvich, D., 2016. Autonomous Identification and Optimal Selection of Popular Smart Vehicles for Urban Sensing - An Information-Centric Approach. IEEE Transactions on Vehicular Technology 65 (12), 9529-9541.

\bibitem[Kullback and Leibler, 1951]{KL1951} Kullback, S., Leibler, R.A., 1951. On information and sufficiency. The annals of mathematical statistics, 22(1), 79-86.

\bibitem[Kumar et al., 2015]{Kumar2015} Kumar, P., Morawska, L., Martani, C., Biskos, G., Neophytou, M., Di Sabatino, S., Bell, M., Norford, L., Britter, R., 2015. The rise of low-cost sensing for managing air pollution in cities. Environment International. 75, 199-205.

\bibitem[Li et al., 2009]{LSLHLW2009} Li, X., Shu, W., Li, M., Huang,H.Y., Luo, P.E., Wu,M.Y.,2009. Performance Evaluation of Vehicle-Based Mobile Sensor Networks for Traffic Monitoring. IEEE Transactions on Vehicular Technology, 58(4), 1647-1653. 

\bibitem[Li et al., 2023]{LMZLL2023} Li Y., Meng X., Zhao H., Li W., Long Y., 2023. Identifying abandoned buildings in shrinking cities with mobile sensing images. Urban Informatics, 2, 3.

\bibitem[Liu et al., 2005]{LBDNT2005} Liu, B., Brass, P., Dousse, O., Nain, P., Towsley, D., 2005. Mobility improves coverage of sensor networks. in: Proceedings of the 6th ACM international symposium on Mobile ad hoc networking and computing  - MobiHoc ’05. New York, NY, USA, 300-308. 

\bibitem[Liu et al., 2016]{LNL2016}  Liu, Y., Niu, J., Liu, X., 2016. Comprehensive tempo-spatial data collection in crowd sensing using a heterogeneous sensing vehicle selection method. Personal and Ubiquitous Computing 20 (3), 397-411.

\bibitem[Long et al., 2019]{LSPHZL2019} Long, J., Sun, Z., Pardalos, P.M., Hong, Y., Zhang, S., Li, C., 2019. A hybrid multi-objective genetic local search algorithm for the prize-collecting vehicle routing problem. Information Sciences 478, 40-61.


\bibitem[Ma et al., 2014]{MZY2014} Ma, H., Zhao, D., Yuan, P., 2014. Opportunities in mobile crowd sensing. IEEE Communications Magazine 52 (8), 29-35.

\bibitem[Maag et al., 2017]{Maag2017} Maag, B., Zhou, Z., Saukh, O., Thiele, L., 2017. SCAN: Multi-hop calibration for mobile sensor arrays. in:  Proceedings of the ACM on Interactive, Mobile, Wearable and Ubiquitous Technologies, 19. 

\bibitem[Meegahapola et al., 2019]{Meegahapola2019} Meegahapola, L., Kandappu, T., Jayarajah, K., Akoglu, L., Xiang, S., Misra, A., 2019. BuScope: Fusing Individual \& Aggregated Mobility Behavior for “Live” Smart City Services, in: Proceedings of the 17th Annual International Conference on Mobile Systems, Applications, and Services. Presented at the MobiSys ’19: The 17th Annual International Conference on Mobile Systems, Applications, and Services, ACM, New York, NY, USA, pp. 41–53. 

\bibitem[O'Keeffe et al., 2019]{OASSR2019} O'Keeffe, K.P., Anjomshoaa, A., Strogatz, S.H., Santi, P., Ratti, C., 2019. Quantifying the sensing power of vehicle fleets. Proceedings of the National Academy of Sciences 116 (26), 12752–12757.

\bibitem[Pasqualetti et al., 2012]{PDB2012} Pasqualetti, F., Durham, J.W., Bullo, F., 2012. Cooperative Patrolling via Weighted Tours: Performance Analysis and Distributed Algorithms. IEEE Transactions on Robotics 28 (5), 1181-1188.

\bibitem[Riahi et al., 2021]{RHS2021} Riahi, V., Hakim Newton, M.A., Sattar, A., 2021. A scatter search algorithm for time-dependent prize-collecting arc routing problems. Computers \& Operations Research 134, 105392.

\bibitem[Rifki et al., 2020]{RCS2020} Rifki, O., Chiabaut, N., Solnon, R., 2020. On the impact of spatio-temporal granularity of traffic conditions on the quality of pickup and delivery optimal tours. Transportation Research Part E: Logistics and Transportation Review, 142:102085.


\bibitem[Rodriguez-Vega et al., 2019]{RCF2019} Rodriguez-Vega, M., Canudas-de-Wit, C., Fourati, H., 2019. Location of turning ratio and flow sensors for flow reconstruction in large traffic networks. Transportation Research Part B: Methodological 121, 21-40.


\bibitem[Salari et al., 2019]{SKLLE2019} Salari, M., Kattan, L., Lam, W.H.K., Lo, H.P., Esfeh, M.A., 2019. Optimization of traffic sensor location for complete link flow observability in traffic network considering sensor failure. Transportation Research Part B: Methodological 121, 216-251.

\bibitem[Saukh et al., 2012]{Saukh2012} Saukh, O., Hasenfratz, D., Noori, A., Ulrich, T., Thiele, L., 2012. Route selection for mobile sensors with checkpointing constraints, in: 2012 IEEE International Conference on Pervasive Computing and Communications Workshops. Presented at the 2012 IEEE International Conference on Pervasive Computing and Communications Workshops (PerCom Workshops), IEEE, Lugano, Switzerland, pp. 266–271. 

\bibitem[Saukh et al., 2014]{Saukh2014} Saukh, O., Hasenfratz, D., Walser, C., Thiele, L., 2014. On rendezvous in mobile sensing networks.  Lecture Notes in Electrical Engineering 281, 29-42.

\bibitem[Song et al., 2021]{SHS2021} Song, J., Han, K., Stettler, M.E.J., 2021. Deep-MAPS: Machine learning based mobile air pollution sensing. IEEE Internet of Things Journal 8 (9), 7649-7660.




\bibitem[Tonekaboni et al., 2020]{TRMSO2020} Tonekaboni, N.H., Ramaswamy, L., Mishra, D., Setayeshfar, O., Omidvar, S., 2020. Spatio-Temporal Coverage Enhancement in Drive-By Sensing Through Utility-Aware Mobile Agent Selection. in: Internet of Things - ICIOT 2020, vol. 12405, W. Song, K. Lee, Z. Yan, L.-J. Zhang, and H. Chen, Eds. Cham: Springer International Publishing, 108-124.

\bibitem[Viti et al., 2014]{VRCT2014} Viti, F., Rinaldi, M., Corman, F., Tamp\'ere, C.M.J., 2014. Assessing partial observability in network sensor location problems. Transportation Research Part B: Methodological 70, 65-89.



\bibitem[Wang et al., 2018]{WLQWL2018} Wang, C., Li, C., Qin, C., Wang, W., Li, X., 2018. Maximizing spatial-temporal coverage in mobile crowd-sensing based on public transports with predictable trajectory. International Journal of Distributed Sensor Networks 14 (8), 155014771879535.

\bibitem[Wang et al., 2014]{Wang2014} Wang, M., Birken, R., Shahini Shamsabadi, S., 2014. Framework and implementation of a continuous network-wide health monitoring system for roadways. in: Nondestructive Characterization for Composite Materials, Aerospace Engineering, Civil Infrastructure, and Homeland Security 2014 (International Society for Optics and Photonics, Bellingham, WA, 2014), vol.9063, p.90630H.

\bibitem[Wang et al., 2021]{WWWZQ2021} Wang, W., Wu, S., Wang, S., Zhen, L., Qu, X., 2021. Emergency facility location problems in logistics: Status and perspectives. Transportation research part E: logistics and transportation review, 154:102465.

\bibitem[Xu et al., 2019]{XCPJZN2019} Xu, S., Chen, X., Pi, X., Joe-Wong, C., Zhang, P., Noh, H.Y., 2019. iLOCuS: Incentivizing Vehicle Mobility to Optimize Sensing Distribution in Crowd Sensing. IEEE Transactions on Mobile Computing 19 (8), 1831-1847.

\bibitem[Xu et al., 2016]{XHCC2016} Xu, X., Lo, H.K., Chen, A., Castillo, E., 2016. Robust network sensor location for complete link flow observability under uncertainty, Transportation Research Part B, 88, 1-20. 






\bibitem[Yi et al., 2017]{YDLCG2017} Yi, K., Du, R., Liu, L., Chen, Q., Gao, K., 2017. Fast participant recruitment algorithm for large-scale Vehicle-based Mobile Crowd Sensing. Pervasive and Mobile Computing 38, 188-199.

\bibitem[Yu et al., 2021]{YFLRL2021} Yu, H., Fang, J., Liu, S., Ren, Y., Lu, J., 2021. A node optimization model based on the spatiotemporal characteristics of the road network for urban traffic mobile crowd sensing. Vehicular Communications 31, 100383.

\bibitem[Yu et al., 2012]{YLL2012} Yu, J.J.Q., Li, V.O.K., Lam, A.Y.S., 2012. Sensor deployment for air pollution monitoring using public transportation system. in: 2012 IEEE Congress on Evolutionary Computation. Brisbane, Australia, 1-7.

\bibitem[Zhang et al., 2018]{ZYL2018} Zhang, X., Yang, Z., Liu, Y., 2018. Vehicle-Based Bi-Objective Crowdsourcing. IEEE Transactions on Intelligent Transportation Systems 19, 3420–3428.

\bibitem[Zhao et al., 2014]{ZML2014} Zhao, D., Ma, H., Liu, L., 2014. Energy-efficient opportunistic coverage for people-centric urban sensing. Wireless Networks 20 (6), 1461-1476. 

\bibitem[Zhao et al., 2015]{ZMLL2015} Zhao, D., Ma, H., Liu, L., Li, X.-Y., 2015. Opportunistic coverage for urban vehicular sensing. Computer Communications 60, 71–85, 2015.


\bibitem[Zhu et al., 2018]{ZFM2018} Zhu, N., Fu, C., Ma, S., 2018. Data-driven distributionally robust optimization approach for reliable travel-time-information-gain-oriented traffic sensor location model. Transportation Research Part B: Methodological 113, 91-120.


\bibitem[Zhu et al., 2022]{ZFZM2022}Zhu, N., Fu, C., Zhang, X., Ma, S., 2022. A network sensor location problem for link flow observability and estimation. European Journal of Operational Research 300 (2), 428-448.





\bibitem[Zhu et al., 2013]{ZLZLZ2013} Zhu, Y., Li, Z., Zhu, H., Li, M., Zhang, Q., 2013. A Compressive Sensing Approach to Urban Traffic Estimation with Probe Vehicles. IEEE Trans. on Mobile Comput. 12, 2289-2302. 


\end{thebibliography}
\end{document}